\documentclass{article}
\usepackage{arxiv}
\usepackage{amsfonts}
\usepackage{amssymb}
\usepackage{mathtools}
\usepackage{amsmath}
\usepackage{algpseudocode}
\usepackage{algorithm}
\usepackage{comment}
\usepackage{booktabs}
\usepackage{multirow}
\usepackage{tabularx}
\usepackage{xcolor}
\usepackage{colortbl}
\usepackage{pgfplots}
\usepackage{tikz}
\usepackage{needspace}
\usepackage{enumitem}
\usepgfplotslibrary{groupplots,fillbetween}
\usetikzlibrary{decorations.pathreplacing}
\pgfplotsset{compat=1.18}

  \usepackage{amsthm}
  \usepackage{placeins}
  \usepackage[hidelinks,hypertexnames=false]{hyperref}
  \usepackage[capitalize,nameinlink]{cleveref}

  \newtheorem{theorem}{Theorem}[section]
  \newtheorem{lemma}[theorem]{Lemma}
  
  \theoremstyle{definition}
  \newtheorem{definition}[theorem]{Definition}
  \theoremstyle{remark}
  
  \AddToHook{env/lemma/begin}{\crefalias{theorem}{lemma}}
  \AddToHook{env/proposition/begin}{\crefalias{theorem}{proposition}}
  \AddToHook{env/definition/begin}{\crefalias{theorem}{definition}}
  \AddToHook{env/remark/begin}{\crefalias{theorem}{remark}}
  \numberwithin{equation}{section}
  \numberwithin{figure}{section}
  \numberwithin{table}{section}
  \numberwithin{algorithm}{section}

\newcommand{\F}{\mathbb F}
\newcommand{\cP}{\mathcal P}
\newcommand{\epsmach}{\epsilon_{\mathrm{mach}}}
\newcommand{\kgram}{\textup{\textsc{SyrkBatched}}}
\newcommand{\kevd}{\textup{\textsc{SyevjBatched}}}
\newcommand{\krot}{\textup{\textsc{GemmBatched}}}

\newcommand{\mfam}[1]{\emph{#1}}

\newcolumntype{Y}{>{\raggedright\arraybackslash}X}

\crefname{section}{\S}{\S\S}
\Crefname{section}{\S}{\S\S}
\crefname{subsection}{\S}{\S\S}
\Crefname{subsection}{\S}{\S\S}
\crefname{subsubsection}{\S}{\S\S}
\Crefname{subsubsection}{\S}{\S\S}

\newcommand{\sectionrefformat}[1]{%
  \crefformat{#1}{\S##2##1##3}\Crefformat{#1}{\S##2##1##3}%
  \crefmultiformat{#1}{\S\S##2##1##3}{ and~##2##1##3}{, ##2##1##3}{, and~##2##1##3}%
  \Crefmultiformat{#1}{\S\S##2##1##3}{ and~##2##1##3}{, ##2##1##3}{, and~##2##1##3}}
\AtBeginDocument{%
  \sectionrefformat{section}\sectionrefformat{subsection}%
  \sectionrefformat{subsubsection}%
  }

\newcommand{\papershorttitle}{Parallel Block Jacobi with Variable Partitioning}

\newcommand{\papertitle}{Parallel Block Jacobi SVD with Variable Column Partitioning}
\title{\papertitle}
\AtBeginDocument{\hypersetup{pdftitle={\papertitle},
  pdfauthor={Alessandro Carraro, Francesco Tudisco}}}
\newcommand{\paperdepartment}{School of Mathematics and Maxwell Institute for
  Mathematical Sciences}

\newcommand{\correspondingemail}{a.carraro@sms.ed.ac.uk}
\newcommand{\coauthoremail}{f.tudisco@ed.ac.uk}
\newcommand{\coderepository}{\url{https://github.com/alecarraro/ag-osbj}}

  \author{Alessandro Carraro \qquad Francesco Tudisco\\
    \normalfont\paperdepartment\\
    \normalfont University of Edinburgh\\
    \normalfont\href{mailto:\correspondingemail}{\correspondingemail}
    (corresponding author), \href{mailto:\coauthoremail}{\coauthoremail}}

  \renewcommand{\shorttitle}{\papershorttitle}
  \date{}

\begin{document}
\maketitle

\begin{abstract}
Parallel one-sided Jacobi methods for the singular value
decomposition (SVD) organize column orthogonalizations into sweeps
of parallel stages, with each sweep covering every column pair.
Conventional block Jacobi methods guarantee this coverage by visiting
every pair of blocks in a fixed column partition, at the cost of
repeatedly processing within-block column pairs.
We introduce \emph{variable partitioning}, in which the column groups
change between stages and are still processed as blocks. It repeats
fewer column pairs and gives sweeps with nearly half as many stages as conventional fixed-partition
schedules, attaining the minimum possible stage count in several
parameter ranges.
We implement these schedules in AG-OSBJ, a double-precision GPU solver.
Experiments show that the shorter sweeps generally require a similar
number of sweeps to converge, reducing the total number of stages
to convergence in $1320$ of $1326$ matched runs.
Dynamic block ordering can further reduce the stage count for small
matrices, but its advantage diminishes with matrix size, and it needs
many more inner Jacobi sweeps and scheduling work before every stage.
On an NVIDIA H200 NVL, AG-OSBJ achieves median speedups
of $1.6$ to $2.0$ over fixed partitioning with round-robin scheduling.
It is at least $3.7$ times faster than cuSOLVER's Jacobi routine
\texttt{gesvdj}, faster than \texttt{gesvd} on nearly every test matrix,
and faster than \texttt{gesvdp} on most matrices up to order $2048$.
Its errors are close to roundoff level, and it computes the singular
values of strongly column-graded matrices to high relative accuracy.
\end{abstract}

\section{Introduction}
\label{sec:intro}
One-sided Jacobi methods compute the singular value decomposition (SVD)
of a matrix $A\in\mathbb{R}^{m\times n}$ by progressively orthogonalizing
its columns~\cite{Hestenes1958}; once the columns are orthogonal, their
norms are the singular values. These methods are known for their accuracy:
for column-graded matrices, they compute the singular values to high
relative accuracy~\cite{DemmelVeselic1992,DrmacVeselic2008a,DrmacVeselic2008b}.
Their cost, however, has traditionally limited their use for large matrices.

Parallelism can offset this cost. In the scalar one-sided Jacobi
method, a pass called a \emph{sweep} visits every pair of columns.
Operations on disjoint column sets are independent, so for parallel
execution they are grouped into \emph{stages}: operations within a
stage execute concurrently, whereas successive stages execute
sequentially. Block Jacobi methods improve the
efficiency of each stage by processing groups of columns through
matrix--matrix operations~\cite{BeckaVajtersic1999,Novakovic2015,HariBegovic2017}.
Conventionally, the columns are divided into fixed blocks, and
a sweep visits every pair of blocks. Each subproblem processes the
union of two blocks, whose column indices we call a \emph{pivot set}.
With equal blocks of $K/2$ columns and an even number of blocks,
round-robin scheduling of the block pairs, \emph{block round-robin},
requires $2n/K-1$ stages per sweep, each processing
pivot sets of $K$ columns~\cite{BrentLuk1985,LukPark1989,Novakovic2015}.

Visiting every pair of fixed blocks guarantees that every column
pair is processed, but imposes more structure than column-pair
coverage requires. Columns within the same block share a pivot
set at every stage, whereas columns in different blocks meet only
once per sweep, so a panel solve that orthogonalizes all columns
of the pivot set spends a large part of every stage on repeated
within-block pairs (\cref{sec:scalar-orderings}).

In this work, we propose a \emph{variable partitioning} strategy that
separates the grouping of columns for block computations from the
requirement that every column pair be visited during a sweep.
To form such a sweep, we choose a \emph{schedule}: a sequence of
partitions into disjoint pivot sets of at most $K$ columns, allowing the
groups to change between stages while ensuring that every column pair
shares a pivot set at least once. Constructing short sweeps is therefore
a problem of \emph{covering by partitions}~\cite{Furedi1989,FurediGyarfas1991}.
Variable partitioning can nearly halve the stages per sweep; for
example, with $n=2048$ and $K=64$, our schedule requires $33$ stages,
the minimum possible, compared with $63$ for block round-robin.

Partitioning and ordering are separate choices in block Jacobi methods,
and our approach combines variable partitioning with static ordering
(\cref{tab:axes}). Scalar orderings use a fixed partition into
single columns. Static parallel orderings of column or block pairs
have been studied extensively~\cite{BrentLuk1985,ZhouBrent1997,EberleinPark1990,LukPark1989}, with
convergence established through their relation to serial cyclic
orderings~\cite{ForsytheHenrici1960,LukPark1989proof,ShroffSchreiber1989}.
Matrix-dependent strategies instead select interactions during the
iteration. For example, de Rijk's pivoting~\cite{deRijk1989}, used in
LAPACK's \texttt{DGESVJ}~\cite{DrmacVeselic2008b}, selects columns by
their current norms. Dynamic block orderings use estimated block
couplings to select disjoint block pairs through approximate
maximum-weight matchings~\cite{BeckaOksaVajtersic2002,BeckaOksaVajtersic2015,KudoEtAl2017}.
These block orderings retain fixed column groups and their repeated
interactions,\footnote{The ``variable blocking'' of Kudo et al.~\cite{KudoEtAl2017}
changes the number of blocks relative to the number of processors; the
column blocks remain fixed during the iteration.} while computing weights and matchings adds scheduling work before each stage.
For the two-sided Jacobi method, the clique-based ordering of Ok\v{s}a and Vajter\v{s}ic~\cite{OksaVajtersic2006} also changes the pivot sets, which are unions of several blocks of a fixed partition. It selects these unions during the iteration using an approximate maximum-weight clique partitioning computed with a genetic algorithm, and does not arrange them into sweeps that cover every pair of blocks.

Variable partitioning preserves the computational structure of
established high-performance block Jacobi
implementations~\cite{Novakovic2015}.
These implementations orthogonalize the columns of each pivot set
through a small $K\times K$ problem solved in fast on-chip memory, and
apply the accumulated transformation to the full columns as a
matrix--matrix product, which can exploit double-precision (FP64)
tensor cores on recent GPUs~\cite{CUTLASS}.
With variable partitioning, the columns of a pivot set are gathered by
index instead of being read as one contiguous block, and both steps are
unchanged. The schedule does not depend on the matrix entries, so we
construct it when the solver starts, for the given $n$ and $K$, and each
sweep repeats the same sequence of stages without scheduling work
between them.

Our main contributions are as follows:

\begin{table}[t]
\caption{Parallel Jacobi methods by column partitioning (fixed or variable) and pivot ordering (static, pivoted, or dynamic).}
\label{tab:axes}
\centering
\begin{tabular}{@{}lll@{}}
\toprule
& fixed partitioning & variable partitioning \\
\midrule
static ordering
& \cite{BrentLuk1985,LukPark1989,ShroffSchreiber1989,BeckaVajtersic1999,Novakovic2015}
& this paper \\
pivoted or dynamic ordering
& \cite{deRijk1989},~\cite{BeckaOksaVajtersic2002,BeckaOksaVajtersic2015,KudoEtAl2017}
& \cite{OksaVajtersic2006} \\
\bottomrule
\end{tabular}
\end{table}

\begin{itemize}[leftmargin=*]
\item
We propose variable partitioning for the one-sided block Jacobi SVD
and formulate the construction of short static sweeps as a covering
problem (\cref{sec:orderings}).

\item
Drawing on finite geometry and design theory, we construct schedules
for arbitrary column counts $n$ and maximum pivot-set sizes $K$.
We establish optimality for several parameter ranges, including
powers-of-two $n$ and $K$, bound the excess stage count for
$n\leq K^2$, and show that the relative excess vanishes as
$n$ grows (\cref{sec:construction}).

\item
We show experimentally that fewer stages per sweep generally
translate into fewer stages to convergence, reducing the total
stage count in $1320$ of $1326$ matched runs. Compared with dynamic
block ordering, our method requires less computational work per
stage and achieves comparable stage counts as $n$ grows
(\cref{sec:exp-sweeps,sec:exp-dynamic}).

\item
We develop AG-OSBJ, a double-precision GPU solver that turns these
savings into runtime speedups over fixed partitioning with
round-robin scheduling. It is faster than NVIDIA cuSOLVER's
\texttt{gesvdj} in every test and than \texttt{gesvd} in nearly
every one, is competitive with the faster \texttt{gesvdp}, and
computes the singular values of strongly column-graded matrices to
high relative accuracy
(\cref{sec:implementation,sec:exp-sweeps,sec:exp-cusolver}).
\end{itemize}

\section{Sweeps and lower bounds}
\label{sec:orderings}

In this section, we first review the one-sided block Jacobi method with
a fixed column partition, then formulate variable partitioning as a
covering problem and recall lower bounds on the number of stages
required to cover every pair of columns.

\subsection{Fixed partitioning}
\label{sec:scalar-orderings}

Let \(A\in\mathbb{R}^{m\times n}\), \(m\ge n\), be the input matrix,
and initialize the working matrices \(X=A\) and \(V=I_n\). The columns
\(x_1,\dots,x_n\) of \(X\) are indexed by \([n]=\{1,\dots,n\}\), and
for \(S\subseteq[n]\), \(X_S\) and \(V_S\) denote the columns of \(X\)
and \(V\) indexed by \(S\).
In the one-sided block Jacobi method, a fixed partition
$\mathcal B=\{B_1,\dots,B_\ell\}$ of $[n]$ defines the block columns~\cite{HariBegovic2017}.
A \emph{pivot pair} $\{i,j\}$ of block indices  determines the pivot set
$S=B_i\cup B_j$ and the corresponding \emph{pivot Gram matrix}
\[
    H_S=(X^\top X)_{S,S}=X_S^\top X_S.
\]
Processing the pivot pair consists of computing an orthogonal
transformation $Q_S$ that (approximately) diagonalizes $H_S$,
and updating $X_S\leftarrow X_SQ_S$ and $V_S\leftarrow V_SQ_S$.
Pivot pairs are processed repeatedly, in sweeps defined below, until
\begin{equation}
\label{eq:stopping}
    |x_u^\top x_v|
    \le \mathrm{tol}\,\|x_u\|\,\|x_v\|
    \qquad\text{for all }u\ne v,
\end{equation}
where $\mathrm{tol}$ is a small multiple of the machine epsilon
$\epsmach=2^{-52}$, typically $\mathrm{tol}=m\epsmach$~\cite{DemmelVeselic1992,DrmacVeselic2008a}.
The singular values are then the column norms of $X$.

A cyclic block Jacobi strategy repeats a fixed sequence of pivot pairs
in which every pair $\{i,j\}$, $1\le i<j\le\ell$, is selected exactly
once. One pass through this sequence is a sweep.
For parallel execution, disjoint pivot pairs are processed simultaneously in
parallel steps, also called stages~\cite[\S3]{Novakovic2015}.
We write a parallel sweep as
$\cP_1,\dots,\cP_s$, where $\cP_k$ is the set of pivot sets
$B_i\cup B_j$ processed at stage $k$. \Cref{alg:block-jacobi}
summarizes the one-sided block Jacobi algorithm for a fixed column
partition and a static ordering.

\begin{algorithm}[t]
\caption{One-sided block Jacobi SVD with fixed partition and static ordering.}
\label{alg:block-jacobi}
\begin{algorithmic}[1]
\Require $X\in\mathbb{R}^{m\times n}$, partition $\mathcal B=\{B_1,\dots,B_\ell\}$,
         sweep $\cP_1,\dots,\cP_s$
\State $V\gets I_n$
\Repeat
    \For{$k=1,\dots,s$}
        \ForAll{$S\in\cP_k$ \textbf{in parallel}}
            \State Compute orthogonal $Q_S$ that (approximately)
                   diagonalizes $X_S^\top X_S$
            \State $X_S\gets X_SQ_S$, $V_S\gets V_SQ_S$
        \EndFor
    \EndFor
\Until{\eqref{eq:stopping} holds}
\State $\sigma_j\gets\|x_j\|$ and
       $u_j\gets x_j/\sigma_j$ for $\sigma_j>0$
\State Complete $U$ if necessary and
       \Return $U,\Sigma=\operatorname{diag}(\sigma_j),V$
\end{algorithmic}
\end{algorithm}

For equal block sizes $|B_i|=K/2$ and even $\ell=2n/K$, every block must
be paired with each of the $\ell-1$ others, and with at most one of them
per stage, so a parallel sweep has at least $\ell-1=2n/K-1$ stages.
Many orderings attain this bound; we use block
round-robin~\cite{BrentLuk1985} as the representative in our comparisons, but the discussion below
applies to every fixed-partition sweep. The within-block pairs, which
share a pivot set at every stage, make up a fraction $(K-2)/(2K-2)$ of
the pairs of each pivot set, whatever the ordering; this fraction
approaches $1/2$ as $K$ grows. A panel solve that diagonalizes $H_S$
processes them at every stage. Zhou and Brent~\cite{ZhouBrent1997}
skip these repetitions; each stage then processes only the $K^2/4$
cross-block pairs of each pivot set, but a sweep still needs $2n/K-1$
sequential stages.

\subsection{Variable partitioning}
\label{sec:variable}

We generalize the notions of stages and sweeps by allowing the column
partition to change between stages. Rather than requiring every pair
of fixed blocks to interact, we require only that every pair of
columns share a pivot set at least once per sweep.

\begin{definition}[Stage and sweep]
\label{def:stage}
For $K\ge2$, a \emph{stage} is a partition $\cP$ of $[n]$ into
\emph{pivot sets} $S$ with $|S|\le K$. A \emph{sweep of length $s$} is a
sequence of stages $\cP_1,\dots,\cP_s$ such that for all
$1\le u<v\le n$ there exist $k$ and $S\in\cP_k$ with
$\{u,v\}\subseteq S$.
\end{definition}

\Cref{def:stage} makes the search for short sweeps a pair-covering
problem. Such problems have been studied, in somewhat different
languages, in the literature on covers by partitions~\cite{FurediGyarfas1991},
on edge colorings~\cite{Furedi1989}, and in design theory~\cite{Haemers1999}.
F\"uredi and Gy\'arf\'as~\cite{FurediGyarfas1991} call a sweep of length $s$ a $2$-cover of $[n]$
by $s$ partitions of rank at most $K$. With the columns as the vertices
of the complete graph on $[n]$, coloring each edge with a stage that
covers it turns a sweep into an edge coloring with $s$ colors whose
monochromatic components have at most $K$ vertices. F\"uredi calls these
$(D,c)$-colorings~\cite[\S3]{Furedi1989}, with $D=K$ and $c=s$.
In design theory, a sweep is a resolvable pair covering with blocks of
size at most $K$, where \emph{resolvable} means that the blocks can be
arranged into \emph{parallel classes}, each partitioning all points~\cite{HughesPiper1985}; here the blocks are the pivot sets, not the blocks
$B_i$ of a fixed partition. With blocks of size exactly $K$, these are
called resolvable coverings~\cite{Haemers1999,vanDamHaemersPeek2003},
and if every pair is moreover covered exactly once, they are resolvable
$2$-$(n,K,1)$ designs~\cite{HughesPiper1985}. Gy\'arf\'as~\cite{Gyarfas2010}
surveys covers by partitions and colorings. \Cref{tab:dictionary}
links our terminology to those of covers by partitions and of design theory.

\begin{table}[t]
\caption{Corresponding terms in block Jacobi methods, covers by partitions, and design theory.}
\label{tab:dictionary}
\centering
\begin{tabularx}{\linewidth}{@{}>{\hsize=.72\hsize}Y>{\hsize=.94\hsize}Y>{\hsize=1.34\hsize}Y@{}}
\toprule
Block Jacobi & Covers by partitions & Designs \\
\midrule
column & element & point \\
pivot set & part & block \\
stage & partition & parallel class \\
sweep of length $s$ & $2$-cover by $s$ partitions & resolvable covering with $s$ classes \\
$|S|\le K$ & rank $\le K$ & blocks of size $\le K$ \\
\bottomrule
\end{tabularx}
\end{table}

\begin{figure}[t]
    \centering
    \input{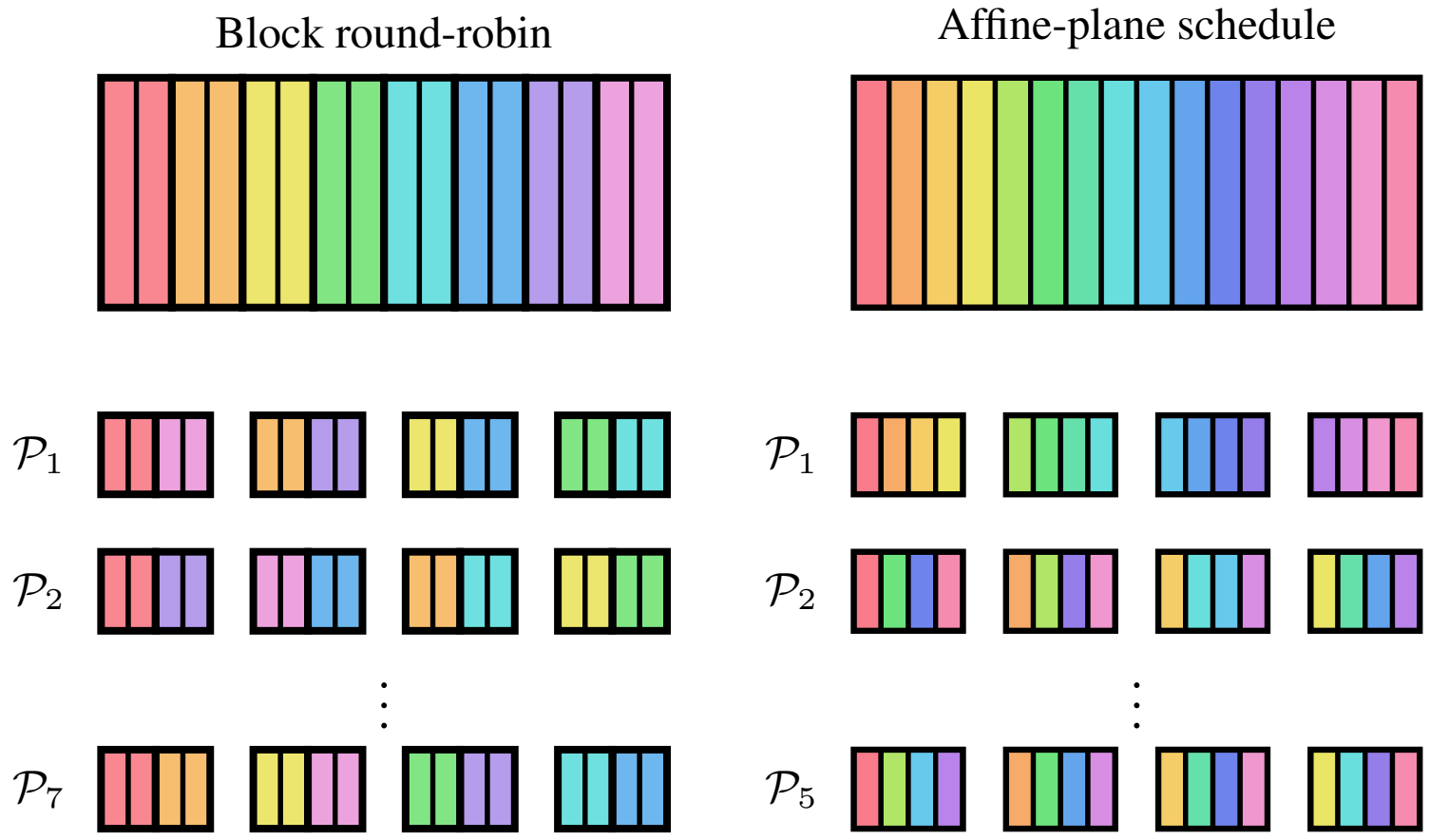}
    \caption{The $n=16$ columns of a matrix (top) and the pivot sets of
    three stages for $K=4$: block round-robin over eight fixed blocks of two
    columns (left, $7$ stages per sweep) and the affine plane
    $\mathrm{AG}(2,4)$ (right, $5$ stages per sweep). Colors identify
    blocks (left) and columns (right).}
    \label{fig:grouping}
\end{figure}

In this language, a shortest fixed-partition sweep with even $\ell$ is a
one-factorization of the complete graph on the blocks $B_1,\dots,B_\ell$,
that is, a resolvable $2$-$(\ell,2,1)$ design, in which every point is
replaced by the $K/2$ columns of its block. \Cref{fig:grouping} compares
such a sweep with a sweep of variable partitioning.

\subsection{Lower bounds}
\label{sec:lower-bounds}
Let $n>K\ge2$, and let $q=\lceil n/K\rceil$ denote the minimum number of pivot sets in a stage. Since each column can meet at most $K-1$ other columns per stage, covering all its $n-1$ partners requires a sweep of length
\begin{equation}
    \label{eq:counting-bound}
    s\ge L(n,K):=\left\lceil\frac{n-1}{K-1}\right\rceil.
\end{equation}
At fixed $K$ and $n\to\infty$, this bound grows as $n/(K-1)$. In comparison,
a fixed partition with even $K$ needs $2q-1$ stages, so their ratio tends to $K/[2(K-1)]$,
which is approximately one half for large enough $K$.

The counting bound ignores that every stage partitions
all the columns. Using this property, Gy\'arf\'as~\cite{Gyarfas1977} shows
(see also F\"uredi~\cite{Furedi1981}) that
\begin{equation}
\label{thm:plane-bound}
    n\le(s-1)K\quad \Longleftrightarrow \quad s\ge q+1.
\end{equation}
For $q\le K$, this exceeds $L(n,K)$ by one exactly when $n\le q(K-1)+1$; for
larger $q$ it adds nothing. Combining the two bounds gives
\begin{equation}
\label{eq:bound}
    s\ge B(n,K):=\max\bigl\{L(n,K),\ q+1\bigr\}
    =\begin{cases}
        q+1, & K<n\le K^2,\\
        L(n,K), & n>K^2.
    \end{cases}
\end{equation}
Sharper bounds are available for specific parameter
regimes~\cite{BierbrauerGyarfas1987,Furedi1989} and may establish the
optimality of schedules that exceed $B(n,K)$. However, the schedules we
construct in \cref{sec:construction} already come close to $B(n,K)$, and a
systematic investigation of sharper bounds is beyond the scope of this
paper.

\section{Constructions}
\label{sec:construction}

Our goal is to construct block Jacobi sweeps that cover every pair of
columns using few sequential stages, with at most $K$ columns in each
pivot set. We begin with schedules based on finite geometry, adapt them
to the required number of columns $n$ and pivot-set size $K$, and select
a sweep for every $n$ and $K$. We then compare the
selected stage count, denoted by $s(n,K)$, with the lower bounds in
\Cref{sec:lower-bounds}. Throughout this section, we maintain the notation
$q=\lceil n/K\rceil$.

\subsection{Affine planes and transversal designs}
\label{sec:geometry}
Our first construction uses the lines of an affine plane as pivot sets
and its parallel classes as stages. Let $\F_p$ be a finite field with
$p$ elements, where $p$ is a prime power. The affine plane
$\mathrm{AG}(2,p)$ has $p^2$ points $(x,y)\in\F_p^2$ and lines
\[
    \lambda=\{y=ax+c\}\quad(a,c\in\F_p),\qquad
    \lambda=\{x=c\}\quad(c\in\F_p).
\]
For each fixed slope $a$, the $p$ lines obtained by varying $c$
partition the points into sets of size $p$. Such a partition into
disjoint lines is called a \emph{parallel class}. The vertical lines
form one more class, giving $p+1$ classes in total.

If we identify the points with matrix columns and use each line as a
pivot set, each parallel class then specifies one stage.
Any two distinct points lie on exactly one common line. Since each line is used as a pivot set in exactly one
stage, every pair of columns shares a pivot set exactly once during
the $p+1$ stages. These stages therefore form a sweep for $n=p^2$ and $K=p$:
in the language of \cref{sec:variable}, a resolvable $2$-$(p^2,p,1)$
design whose blocks are the lines~\cite{HughesPiper1985}.

For the constructions below, it is useful to separate interactions
between groups of columns from those within a group.
A \emph{transversal design} $\mathrm{TD}(v,g)$ partitions $vg$ points
into $v$ groups of size $g$. Its blocks contain one point from each
group and cover every pair from different groups exactly once~\cite{AbelColbournDinitz2007}. We need resolvable ones, whose blocks
form stages.
To obtain such a design from $\mathrm{AG}(2,p)$, we remove the
vertical-line stage and retain its pivot sets as groups.
The remaining $p$ stages cover every pair of points from different
groups and form a resolvable $\mathrm{TD}(p,p)$.
By retaining only $v\le p$ groups and restricting every pivot set to
their points, we obtain a resolvable $\mathrm{TD}(v,p)$.

Only pairs within the same group then remain uncovered.
We cover them by constructing a sweep within each group, a procedure
called \emph{filling in the groups}~\cite{ColbournDinitz2007}.
When these sweeps have the same number of stages, we combine their
$i$th stages into one stage over all columns. This is valid because
the groups are disjoint. We append these stages to those of the
transversal design to complete the sweep.

\subsection{Inflation}
\label{sec:affine}
The plane construction uses pivot sets of size $p$, whereas we
require sizes at most $K$. To adapt it, we assign $w(P)$ matrix columns
to each point $P=(x,y)\in\F_p^2$, where $w(P)$ is a given nonnegative integer and the sets of columns assigned to distinct points are disjoint. We then replace each line by the pivot set
containing all columns assigned to its points. We call this operation
\emph{inflation}, since weights above one enlarge the pivot sets.
Zero weights allow us to remove points. The total number of columns supported,
$\sum_Pw(P)$, is the construction's \emph{capacity}.
F\"uredi uses this point-replacement strategy~\cite[proof of Theorem~3.4]{Furedi1989}, and constant weights recover the
inflated affine plane of Haemers~\cite[Lemma~1]{Haemers1999}.

If we retain all lines, columns assigned to distinct points share
the pivot set associated with the line through those points, while
columns assigned to the same point share the pivot sets associated
with all lines through that point.
Inflation therefore preserves pair coverage, as well as the partition
into disjoint pivot sets at each stage.
As in \cref{sec:geometry}, we use the nonvertical lines as pivot sets
and reserve the vertical lines as groups whose internal interactions
we will schedule separately. We therefore impose the size limit $K$
on each nonvertical line by requiring
\begin{equation}
\label{eq:affine-load}
    \sum_{P\in\lambda}w(P)\le K
    \qquad\text{for every nonvertical line }\lambda.
\end{equation}
The left-hand side counts the columns in the pivot set associated
with line $\lambda$. A vertical line $\lambda$ instead defines a group
of $\sum_{P\in\lambda}w(P)$ columns. Under \eqref{eq:affine-load}, the
nonvertical parallel classes specify $p$ stages satisfying the size
limit, with empty pivot sets omitted.
These stages cover every pair from different groups, and we complete
coverage within groups in one of two ways (\cref{tab:constructions}).

\begin{itemize}[leftmargin=*]
\item \emph{One filling stage.}
If each group has at most $K$ columns, we use the groups themselves
as the pivot sets of one additional stage, so the sweep requires $p+1$
stages in all.
Every line then carries at most $K$ columns; a weight assignment with
this property is called a \emph{weighted $K$-arc}~\cite{Bierbrauer1991}.
If $p$ divides $K$, the constant weights $w(P)=K/p$ put $K$ columns on
every line and support $pK$ columns. For every other prime power $p<K$,
the weights of \cref{app:weights} support $C(p,K)=p(K-1)$ columns.

\item \emph{Recursive filling.}
For $p\ge K$, we assign one column to each point of $K$ vertical lines
and none elsewhere. Each nonvertical line meets each of these lines once,
so its pivot set has $K$ columns, and the $p$ stages together support
$pK$ columns. We number the columns along the horizontal lines, which
leaves at most $q$ columns in each group. We then construct a sweep on
$q$ columns by the same procedure and apply a copy within each group,
combining corresponding stages across groups as in \cref{sec:geometry}.
Appending these $s(q,K)$ stages to the $p$ stages between groups gives a
sweep of length $p+s(q,K)$, where $s(q,K)$ is the number of stages in the
sweep selected for $q$ columns.
\end{itemize}

\begin{table}[t]
\caption{Sweeps from an affine plane of order $p$ for pivot sets of at
most $K$ columns: assignment of columns to points, condition on $p$,
stages per sweep, and the largest number of columns covered.}
\label{tab:constructions}
\centering
\small
\begin{tabularx}{\linewidth}{@{}Ylcc@{}}
\toprule
Weights & Order & Stages & Capacity \\
\midrule
$K/p$ columns on every point & $p\mid K$ & $p+1$ & $pK$ \\
Weighted $K$-arc (\cref{app:weights}) & $p\nmid K$, $p<K$ & $p+1$ & $p(K-1)$ \\
One column on the points of $K$ vertical lines & $p\ge K$ & $p+s(q,K)$ & $pK$ \\
\bottomrule
\end{tabularx}
\end{table}

\subsection{Selection}
\label{sec:choice}

To select a sweep for $n$ columns, we take the smallest prime-power order
$p\ge q$ in \cref{tab:constructions} whose capacity is at least $n$.
A larger capacity suffices because deleting excess columns from every
pivot set preserves coverage of the remaining pairs.

For $n\le K$, we use one pivot set. For $K<n\le K^2$, one filling stage
gives $p+1$ stages, where $p\le K$ is the smallest prime power with
$p\ge q$ that holds $n$ columns: $n\le pK$ if $p$ divides $K$, and
$n\le C(p,K)$ otherwise. When $K$ is a power of two, as in our
implementation, such a $p$ always exists, since the plane of order $K$
holds $K^2$ columns. For other $K$, if there is none, recursive filling
with the smallest prime power $p\ge K$ also gives $p+1$ stages, because
its groups have at most $q\le K$ columns.

For $n>K^2$, we use recursive filling with the smallest prime power
$p\ge q$, so that
\begin{equation}
\label{eq:large-recursion}
    s(n,K)=p+s(q,K).
\end{equation}
Since $q<n$, the recursion reaches the case $n\le K^2$ and defines
$s(n,K)$ for every $n$.

Although these sweeps have few stages, each stage may require
more arithmetic work. For example, when $p>q$, a stage may contain up to
$p$ pivot sets, compared with $q$ for a fixed partition.
In the implementation, we gather the columns of each pivot set into a
submatrix called a \emph{panel} and pad it to $K$ columns, so processing
more pivot sets increases the work even when some contain fewer columns.
Because the pivot sets are disjoint, we can process these panels
concurrently, whereas successive stages must run sequentially.
Our selection rule therefore prioritizes fewer sequential stages while
allowing additional work within each stage.

We measure the excess stage count by $s(n,K)-B(n,K)$, using the lower
bound in \eqref{eq:bound}. For $K<n\le K^2$, the selected sweep has $p+1$
stages and $B(n,K)=q+1$, so the excess is $p-q$. It is zero whenever $q$
is a prime power whose plane holds $n$ columns, that is, $q\mid K$ or
$n\le C(q,K)$, and otherwise at most about $q/5$ for
$q\ge25$~\cite{Nagura1952}.

Above $K^2$, each recursive level may require rounding up to a prime
power. However, prime gaps grow sublinearly~\cite{BakerHarmanPintz2001},
so the excess relative to $B(n,K)$ vanishes as $n\to\infty$ for fixed
$K$.

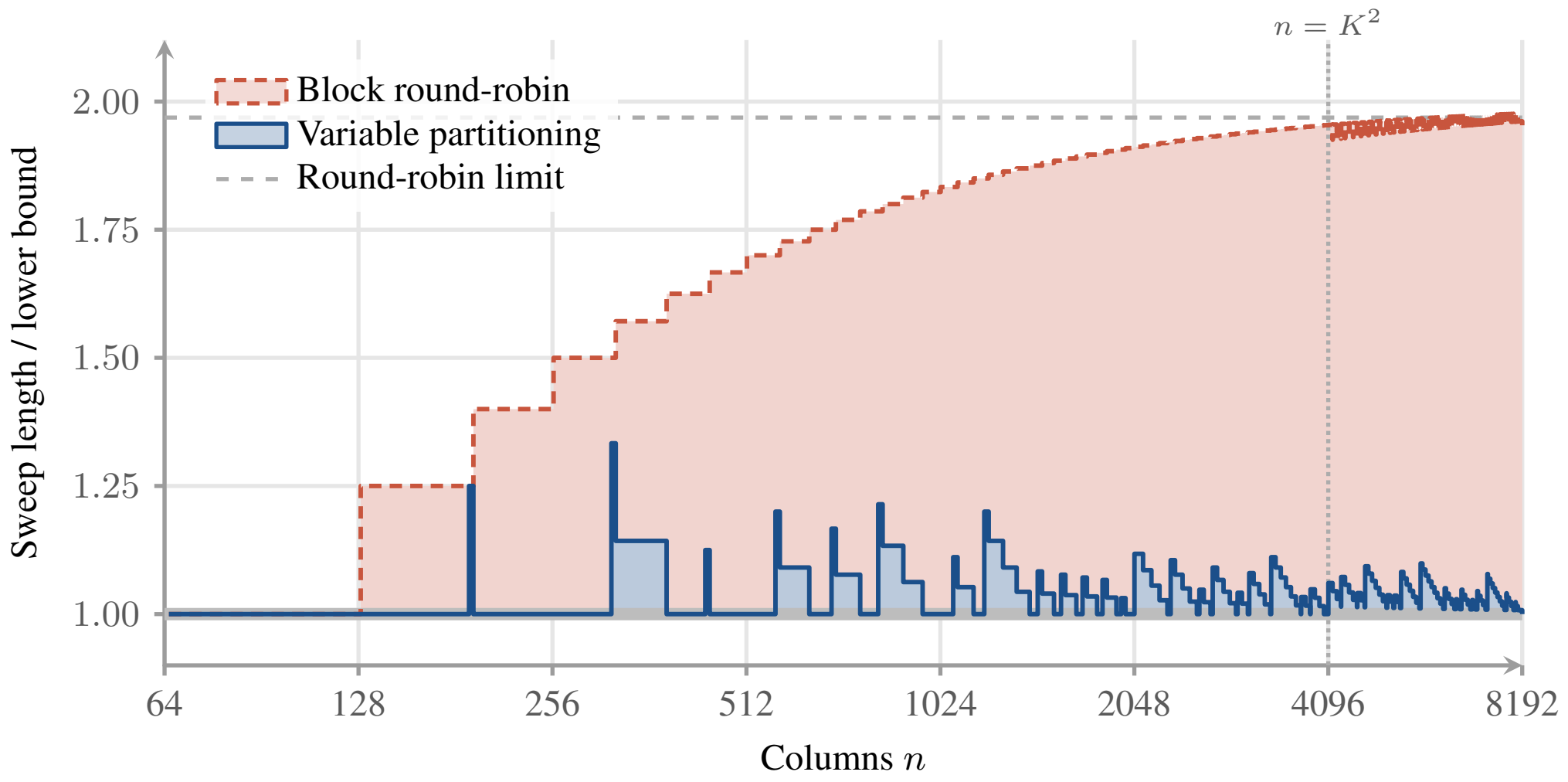
\begin{figure}[t]
\centering
\definecolor{ours}{HTML}{1B4F8A}
\definecolor{robin}{HTML}{C8553D}
\definecolor{bound}{HTML}{5A5A5A}
\begin{tikzpicture}
\begin{axis}[
    width=0.8\linewidth, height=0.42\linewidth,
    axis lines=left, axis line style={bound!60, line width=1bp, line cap=round},
    tick style={bound!60, line width=1bp}, tick align=outside, major tick length=2.5pt,
    ticklabel style={font=\footnotesize, text=bound},
    label style={font=\small},
    grid=major, major grid style={bound!15, line width=1bp},
    xmode=log, log basis x=2,
    xmin=64, xmax=8192,
    xtick={64,128,256,512,1024,2048,4096,8192}, xticklabels={64,128,256,512,1024,2048,4096,8192},
    xlabel={Columns $n$},
    ymin=0.9, ymax=2.12, ytick={1,1.25,1.5,1.75,2},
    yticklabel style={/pgf/number format/.cd, fixed, fixed zerofill, precision=2},
    ylabel={Sweep length / lower bound},
    legend pos=north west,
    legend cell align=left,
    legend style={draw=none, fill=white, fill opacity=0.85, text opacity=1,
        font=\small, row sep=-1pt, inner sep=2pt},
    every axis plot/.append style={line join=round},
    clip=false,
]
\addplot[name path=base, draw=none, forget plot] coordinates {(64,1) (8192,1)};
\addplot[name path=robin, draw=none, forget plot] table[x=n, y=robin] {figures/stages_relative_k64.dat};
\addplot[name path=ours, draw=none, forget plot] table[x=n, y=ours] {figures/stages_relative_k64.dat};
\addplot[fill=robin!25, draw=none, forget plot] fill between[of=robin and ours];
\addplot[fill=ours!30, draw=none, forget plot] fill between[of=ours and base];
\addplot[bound!40, line width=3pt, forget plot] coordinates {(64,1) (8192,1)};
\addplot[bound!50, dashed, line width=1bp, forget plot]
    coordinates {(64,1.96875) (8192,1.96875)};
\addplot[bound!50, densely dotted, line width=1bp, forget plot]
    coordinates {(4096,0.9) (4096,2.12)};
\node[font=\scriptsize, text=bound, anchor=south, inner sep=1pt]
    at (axis cs:4096,2.12) {$n=K^2$};
\addplot[robin, densely dashed, line width=1.1bp, forget plot] table[x=n, y=robin] {figures/stages_relative_k64.dat};
\addplot[ours, line width=1.1bp, forget plot] table[x=n, y=ours] {figures/stages_relative_k64.dat};
\addlegendimage{area legend, line width=1bp, densely dashed, fill=robin!25, draw=robin}
\addlegendentry{Block round-robin}
\addlegendimage{area legend, line width=1bp, fill=ours!30, draw=ours}
\addlegendentry{Variable partitioning}
\addlegendimage{bound!50, dashed, line width=1bp}
\addlegendentry{Round-robin limit}
\end{axis}
\end{tikzpicture}
\caption{Stages per sweep relative to the lower bound $B(n,K)$ against
the number of columns $n$, for variable partitioning and block
round-robin at $K=64$.}
\label{fig:stages-relative}
\end{figure}

When $K$ and $n$ are powers of two, as in most of our
experiments, the excess is zero: below $K^2$ because $q$ divides $K$, and
above because $p=q$ and $L(qK,K)=q+L(q,K)$.
\Cref{fig:stages-relative} compares the stage counts with $B(n,K)$ for
$K=64$ and every $n$ up to $8192$. On average over these sizes, the sweep
is $4\%$ longer than the bound and block round-robin $90\%$ longer.

\section{Implementation}
\label{sec:implementation}

\subsection{Overview}
\label{sec:impl-overview}

Our solver, AG-OSBJ (Affine-Geometry One-Sided Block Jacobi), runs on
NVIDIA Hopper GPUs in FP64. The code is available online.\footnote{
\coderepository}
Each stage consists of three successive kernel launches, each
processing all pivot sets in parallel.
\kgram{} forms $H_S=X_S^\top X_S$, \kevd{} computes a panel rotation
$Q_S$ that approximately diagonalizes it, and \krot{} applies
$[X_S;V_S]\gets[X_S;V_S]\,Q_S$.
We implement two panel widths, $K=64$ and $K=128$, selected according
to the matrix size as explained in \cref{sec:panel-width}.

Because the schedule depends only on $n$ and $K$, the host constructs
it at initialization for one sweep and transfers it to the GPU,
and the sweep loop is captured once in a CUDA graph; together, these
steps take at most about half a sweep.
The iteration terminates when the stopping criterion
\eqref{eq:stopping} is satisfied, with
$\mathrm{tol}=m\epsmach/(1-m\epsmach)$, or when the prescribed
sweep limit is reached. Within a panel, \kevd{} does not rotate a column
pair whose cosine is at most $\sqrt m\,\epsmach$. To reduce its overhead, the convergence check is overlapped with
the beginning of the next sweep, resulting in one additional stage
after convergence is detected.

Two choices shape the solver. First, we omit preconditioning.
Drma\v{c} and Veseli\'c~\cite{DrmacVeselic2008a,DrmacVeselic2008b}
precede the Jacobi iteration with a column-pivoted QR factorization,
which reduces the number of sweeps. However, no efficient GPU
implementation of this factorization is
known~\cite[\S9.2]{AbdelfattahFasi2026}, and therefore GPU Jacobi
solvers typically omit it.
Second, following Abdelfattah and Fasi~\cite{AbdelfattahFasi2026},
we solve the panel problems only approximately: \kevd{} performs
a single Jacobi sweep on $H_S$, rather than diagonalizing it to
convergence. We call this an \emph{inner} sweep, to distinguish
it from the \emph{outer} sweeps of the schedule.

The inner sweep applies variable partitioning within each panel, at a
smaller scale. \kevd{} assigns the $K$ columns to the $64$ points of
$\mathrm{AG}(2,8)$, with one column per point at $K=64$ and two
at $K=128$, and traverses the plane's nine parallel
classes, processing disjoint groups of columns through local
Jacobi orthogonalizations. An implementation overview is provided
in \cref{sm:evd}, with detailed kernel design documentation
available in the GitHub repository.

Since every pair of points lies on a line of some class, the inner
sweep processes every column pair of the panel in exact arithmetic,
so AG-OSBJ is a quasi-cyclic element-wise
Jacobi method~\cite{Hari2007}: its sequence of plane rotations is
periodic, with a period of one outer sweep, and every column pair
occurs at least once per period. With simple singular values, such
methods converge quadratically once the columns are nearly
orthogonal~\cite{Wilkinson1962,Hari2007}.
Existing global convergence analyses of block Jacobi methods with
approximate panel solves assume fixed partitions, restricted partition
changes, or additional conditions on the panel
rotations~\cite{Drmac2009,HariBegovic2017,Hari2025,OksaYamamotoVajtersic2017,OksaYamamotoBeckaVajtersic2018,OksaYamamotoBeckaVajtersic2019},
and do not cover AG-OSBJ. We therefore assess its convergence
empirically in \cref{sec:exp-sweeps}.

In total, each outer sweep applies approximately $n$ rotations to every
column. In exact arithmetic, these rotations keep the columns of $V$ at
unit norm; in floating point, rounding errors accumulate and the norms
drift away from one, which increases the backward error as $n$ grows.
AG-OSBJ reduces this \emph{norm drift} by using Rutishauser's form for
the plane rotations in \kevd{}~\cite{Rutishauser1966}, and removes its
first-order effect at termination by normalizing the columns of $V$ and
scaling those of $X$ by the same factors.
Both measures are analyzed and measured in \cref{sm:drift}.

\subsection{Panel width}
\label{sec:panel-width}
As shown by the numerical experiments in
\cref{sec:exp-sweeps,sec:exp-cusolver}, the number of sweeps required
for convergence and the observed numerical accuracy are largely
insensitive to the choice of panel width $K$.
We therefore choose $K$ for speed alone, to minimize the time of one
outer sweep on the target hardware.

An outer sweep consists of approximately $n/K$ sequential stages,
each containing $n/K$ independent panels.
In each stage, \krot{} performs $4n^2K$ flops on $32n^2$ bytes,
while \kgram{} exploits symmetry to perform $n^2K$ flops on
$8n^2$ bytes.
Both kernels therefore have an arithmetic intensity of $K/8$
flop/byte.
Over a sweep, they perform approximately $5n^3$ flops,
independent of $K$, and move $40n^3/K$ bytes.

On the H200 NVL, the GPU used in our experiments, the peak FP64
tensor-core throughput is $\pi=67\,\mathrm{Tflop/s}$ and the
memory bandwidth is $\beta=4.8\,\mathrm{TB/s}$, placing the
roofline ridge at $K\approx112$.
Increasing $K$ therefore reduces the bandwidth cost until the
compute bound is reached. The cost of \kevd{}, by contrast, increases with $K$.
For $n\leq8192$, all panels within a stage of \kevd{} can be processed concurrently,
requiring at most $128$ of the GPU's $132$ streaming multiprocessors (one per panel
for $K=64$ and two for $K=128$). Consequently, the measured execution time per
stage, $t_{\mathrm E}(K)$, remains approximately independent of $n$, at $45\,\mu\mathrm{s}$ for
$K=64$ and $141\,\mu\mathrm{s}$ for $K=128$. Combining this cost with the roofline model
yields the following asymptotic lower-bound estimate for the sweep time:
\[
T(K)\gtrsim
\max\left\{
    \frac{40n^3}{\beta K},\,
    \frac{5n^3}{\pi}
\right\}
+\frac{n}{K}\,t_{\mathrm E}(K).
\]

For large matrices, the model favors $K=128$, the smallest
supported width above the roofline ridge.
For small matrices, where the working set fits in L2 cache and
kernel latency becomes more important, the lower cost of \kevd{}
at $K=64$ can outweigh the reduction in stage count at $K=128$.
Measurements of the time per sweep place the crossover at $n=3072$,
so AG-OSBJ uses $K=64$ for $n<3072$ and $K=128$ otherwise.

\section{Numerical experiments}
\label{sec:experiments}
\renewcommand\bottomfraction{.6}

In this section, we evaluate the convergence, numerical accuracy, and
running time of AG-OSBJ. We first examine whether shorter sweeps reduce
the stages to convergence without increasing the number of sweeps
(\cref{sec:exp-sweeps}). We then compare with a dynamic ordering
(\cref{sec:exp-dynamic}) and assess the complete solver against
cuSOLVER (\cref{sec:exp-cusolver}).

\subsection{Setup}
\label{sec:exp-setup}

Since the SVD of a rectangular matrix can always be reduced to the square
case by a preliminary QR or LQ factorization, we restrict our
experiments to square matrices. All test matrices have full rank
and belong to the families listed in \cref{tab:matrix-families}.

The \texttt{randsvd} matrices~\cite{Higham2002} are $A=U\Sigma V^\top$,
where $U$ and $V$ are independent Haar-distributed orthogonal matrices,
generated by signed Gaussian QR factorizations, and $\Sigma$ holds one of
the eight spectra in \cref{tab:matrix-families}. Five of these spectra
depend on the condition number $\kappa$, which takes two values, so there
are $13$ \texttt{randsvd} matrices per size. To test relative accuracy, we
also use column-graded matrices $A=BD$. Here $B$ is the \mfam{geometric}
\texttt{randsvd} matrix with $\kappa=10$ and normalized columns, for which
we measure $10\le\kappa_2(B)\le20$, and $D$ is diagonal, with entries from
one of two families with condition number $\kappa_D$, in random order.
This gives four column-graded matrices per size.

\begin{table}[t]
\centering
\small
\caption{Singular values $\sigma_i$ of the \texttt{randsvd} test matrices and
diagonal entries $d_i$ of $D$ in the column-graded test matrices $A=BD$, with
$\theta_i=(i-1)/(n-1)$. The parameter $\kappa$ is the condition number of the
\texttt{randsvd} matrix, $\kappa_D$ that of $D$, and $\mathcal K=\{10^4,10^{10}\}$.
\mfam{Polynomial decay}, \mfam{two plateaus}, \mfam{multi-cluster}, \mfam{one large} and \mfam{one small} follow
\cite{SchwarzEtAl2026}; the $\kappa$-dependent families are the
modes of \texttt{randsvd} and of LAPACK's \texttt{xLATMS}.}
\label{tab:matrix-families}
\renewcommand{\arraystretch}{1.15}
\begin{tabularx}{\linewidth}{@{}llXl@{}}
\toprule
Class & Family & Spectrum or diagonal scaling & Parameter \\
\midrule
\multirow{8}{*}{\texttt{randsvd}}
& \mfam{Polynomial decay}
& $\sigma_i=i^{-2}$
& -- \\

& \mfam{Two plateaus}
& $\sigma_i=10^{-4}+(1+e^{i-10})^{-1}$
& -- \\

& \mfam{Arithmetic}
& $\sigma_i=1-(1-\kappa^{-1})\theta_i$
& $\kappa\in\mathcal K$ \\

& \mfam{Log-random}
& $\sigma_1=1$, $\sigma_n=\kappa^{-1}$;
  $\sigma_{2:n-1}\sim\operatorname{LogUnif}[\kappa^{-1},1]$
& $\kappa\in\mathcal K$ \\

& \mfam{Geometric}
& $\sigma_i=\kappa^{-\theta_i}$
& $\kappa\in\mathcal K$ \\

& \mfam{Multi-cluster}
& $\displaystyle
  \sigma_i=
  \begin{cases}
  n, & 1\le i\le 3,\\
  1, & n-2\le i\le n,\\
  n/2, & n/2-1\le i\le n/2+3,\\
  n+1-i, & \text{otherwise}
  \end{cases}$
& -- \\

& \mfam{One large}
& $\sigma_1=1$, $\sigma_{2:n}=\kappa^{-1}$
& $\kappa\in\mathcal K$ \\

& \mfam{One small}
& $\sigma_{1:n-1}=1$, $\sigma_n=\kappa^{-1}$
& $\kappa\in\mathcal K$ \\

\midrule
\multirow{2}{*}{\begin{tabular}[c]{@{}l@{}}Column-\\graded\end{tabular}}
& \mfam{Geometric}
& $d_i=\kappa_D^{-\theta_i}$
& $\kappa_D\in\mathcal K$ \\

& \mfam{Log-random}
& $d_1=1$, $d_n=\kappa_D^{-1}$;
  $d_{2:n-1}\sim\operatorname{LogUnif}[\kappa_D^{-1},1]$
& $\kappa_D\in\mathcal K$ \\
\bottomrule
\end{tabularx}
\end{table}

The GPU experiments run on an NVIDIA H200 NVL with CUDA 12.9. The GPU
runs at its power limit for about $70\%$ of the time, and its clock varies
between $1215$ and $1785$ MHz. We therefore interleave the solvers so that
all of them run under comparable conditions. For each matrix, we run every
solver once in turn, six times over, with a different solver first each
time, and each run starts from a flushed L2 cache. We discard the first
run of each solver as a warm-up and report the median of the other five;
the slowest of these is within $3\%$ of the fastest in $96\%$ of cases.

\subsection{Sweep length and convergence}
\label{sec:exp-sweeps}

We examine how sweep length affects convergence, using RR-OSBJ as
the baseline. RR-OSBJ (Round-Robin OSBJ) is the same solver with a
fixed partition into blocks of $K/2$ columns and block round-robin
ordering. We denote by $s$ the number of stages per sweep, and by
$N_{\mathrm{stages}}$ and $N_{\mathrm{sweeps}}=N_{\mathrm{stages}}/s$ the
numbers of stages and sweeps to convergence. In this section, we check
the stopping criterion \eqref{eq:stopping} after every stage, so
$N_{\mathrm{sweeps}}$ may be fractional. We use the $13$
\texttt{randsvd} matrices with three seeds each, $39$ test matrices
in all, at both supported panel widths, $K=64$ and $K=128$.
Every run in this section and the next converges, and AG-OSBJ and
RR-OSBJ reach errors of the size reported in \cref{sec:exp-cusolver}.

\begin{figure}[htbp]
\centering
\definecolor{ours}{HTML}{1B4F8A}
\definecolor{robin}{HTML}{C8553D}
\definecolor{bound}{HTML}{5A5A5A}
\definecolor{gesvdj}{HTML}{3A7D44}
\definecolor{gesvdp}{HTML}{C99A2E}
\begin{tikzpicture}
\begin{axis}[
    width=0.68\linewidth, height=0.47\linewidth,
    scale only axis=false,
    axis lines=left, axis line style={bound!60, line width=1bp, line cap=round},
    tick style={bound!60, line width=1bp}, tick align=outside,
    major tick length=2.5pt,
    ticklabel style={font=\footnotesize, text=bound},
    /pgf/number format/1000 sep={},
    label style={font=\small},
    every axis title/.style={font=\small, at={(0.5,1.04)}, anchor=south},
    grid=major, major grid style={bound!15, line width=1bp},
    every axis plot/.append style={line join=round},
    legend cell align=left,
    legend style={draw=none, fill=none,
        font=\scriptsize, row sep=-2pt, inner sep=1pt},
    xmin=0.489, xmax=1.055,
    ymin=0.40, ymax=1.12,
    xtick={0.5,0.6,0.7,0.8,0.9,1},
    xlabel={Relative stages per sweep},
    ylabel={Relative stages to convergence},
    legend style={at={(1.03,0.55)}, anchor=west, font=\footnotesize, row sep=2pt},
]
\addplot[bound, dashed, line width=1bp, domain=0.489:1.055, samples=2] {x};
\addlegendentry{Equal sweep counts}
\addplot[bound!65, densely dotted, line width=1bp, domain=0.489:1.055, samples=2] {1};
\addlegendentry{Equal stage counts}
\addlegendimage{area legend, line width=1bp, draw=ours, fill=ours!20}
\addlegendentry{$K=64$}
\path[draw=ours, fill=ours!20, line width=1bp] (axis cs:0.520310,0.463486) rectangle (axis cs:0.527310,0.496714);
\draw[ours, line width=1.2bp] (axis cs:0.520310,0.480088) -- (axis cs:0.527310,0.480088);
\path[draw=ours, fill=ours!20, line width=1bp] (axis cs:0.552056,0.491501) rectangle (axis cs:0.559056,0.561983);
\draw[ours, line width=1.2bp] (axis cs:0.552056,0.513131) -- (axis cs:0.559056,0.513131);
\path[draw=ours, fill=ours!20, line width=1bp] (axis cs:0.583802,0.515353) rectangle (axis cs:0.590802,0.579133);
\draw[ours, line width=1.2bp] (axis cs:0.583802,0.547312) -- (axis cs:0.590802,0.547312);
\path[draw=ours, fill=ours!20, line width=1bp] (axis cs:0.647294,0.583187) rectangle (axis cs:0.654294,0.654545);
\draw[ours, line width=1.2bp] (axis cs:0.647294,0.608602) -- (axis cs:0.654294,0.608602);
\path[draw=ours, fill=ours!20, line width=1bp] (axis cs:0.774278,0.692809) rectangle (axis cs:0.781278,0.754545);
\draw[ours, line width=1.2bp] (axis cs:0.774278,0.726138) -- (axis cs:0.781278,0.726138);
\addlegendimage{area legend, line width=1bp, draw=robin, fill=robin!20}
\addlegendentry{$K=128$}
\path[draw=robin, fill=robin!20, line width=1bp] (axis cs:0.544887,0.483090) rectangle (axis cs:0.551887,0.538612);
\draw[robin, line width=1.2bp] (axis cs:0.544887,0.504329) -- (axis cs:0.551887,0.504329);
\path[draw=robin, fill=robin!20, line width=1bp] (axis cs:0.609403,0.536105) rectangle (axis cs:0.616403,0.586520);
\draw[robin, line width=1.2bp] (axis cs:0.609403,0.552486) -- (axis cs:0.616403,0.552486);
\path[draw=robin, fill=robin!20, line width=1bp] (axis cs:0.673919,0.592384) rectangle (axis cs:0.680919,0.650644);
\draw[robin, line width=1.2bp] (axis cs:0.673919,0.607735) -- (axis cs:0.680919,0.607735);
\path[draw=robin, fill=robin!20, line width=1bp] (axis cs:0.802952,0.713665) rectangle (axis cs:0.809952,0.777777);
\draw[robin, line width=1.2bp] (axis cs:0.802952,0.747845) -- (axis cs:0.809952,0.747845);
\addplot[only marks, black, mark=square*, mark size=1.5pt, line width=1bp, forget plot] coordinates {(1,1)};
\node[anchor=south east, font=\scriptsize] at (axis cs:0.995,1.015) {RR-OSBJ};
\end{axis}
\end{tikzpicture}
\caption{Stages to convergence against stages per sweep, both relative
to RR-OSBJ, at $n=2048$. One box per sweep length shorter than
RR-OSBJ's, AG-OSBJ's leftmost; RR-OSBJ is the square at $(1,1)$. Each
box gives the median and interquartile range over the $39$ test
matrices ($13$ \texttt{randsvd} matrices, three seeds each).}
\label{fig:sweep-length}
\end{figure}
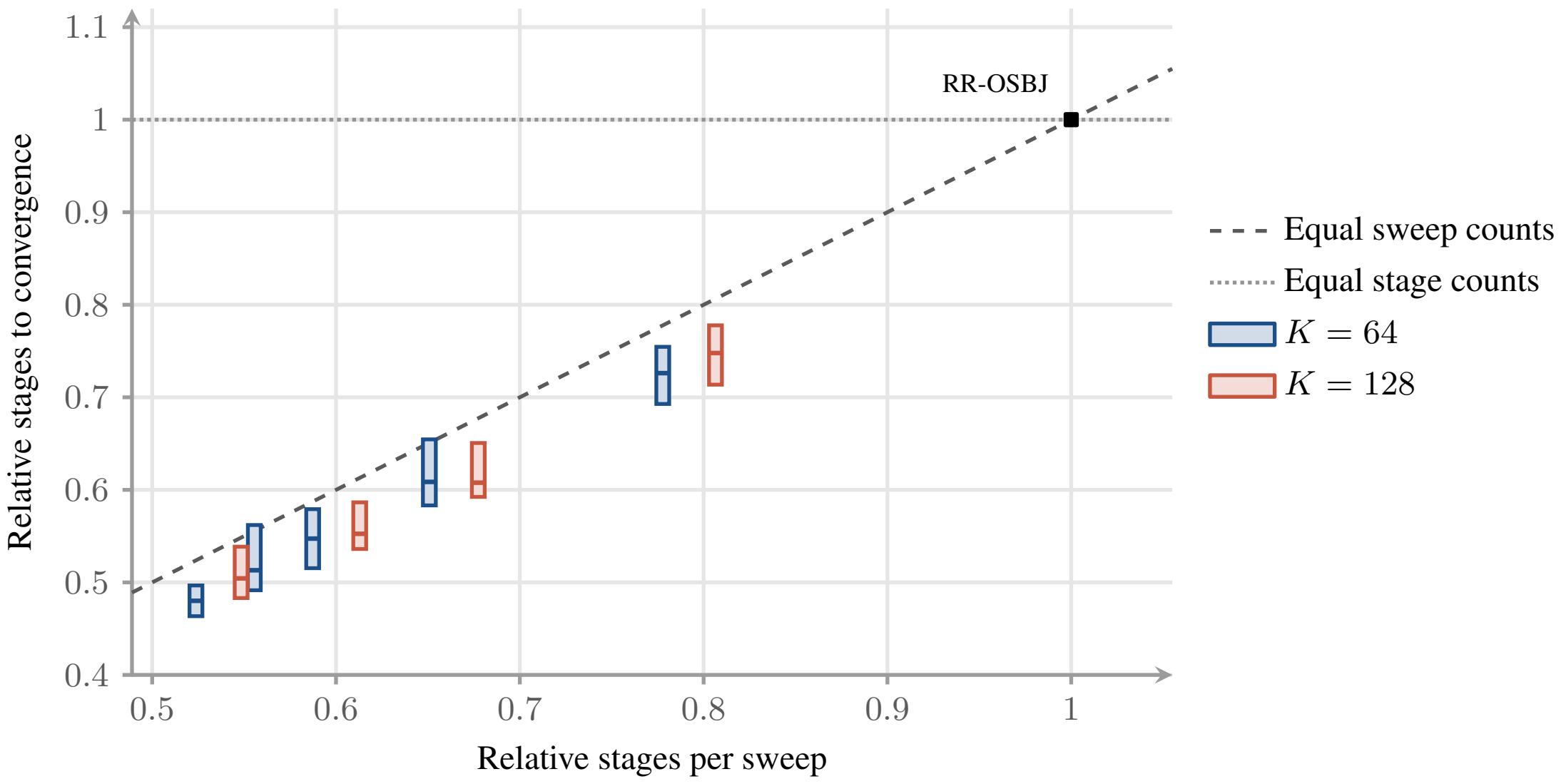

We first vary the sweep length at $n=2048$, keeping the coverage of
every column pair. AG-OSBJ uses $33$ stages per sweep at $K=64$ and
$17$ at $K=128$, and RR-OSBJ uses $63$ and $31$. Recursive filling
(\cref{sec:affine}) with fewer and larger groups gives the
intermediate lengths $35$, $37$, $41$, and $49$ at $K=64$, and $19$,
$21$, and $25$ at $K=128$.
\Cref{fig:sweep-length} compares $N_{\mathrm{stages}}$ against $s$, both normalized by
RR-OSBJ for each matched matrix and seed.

The main observation is that the total number of stages to
convergence is nearly proportional to the number of stages per sweep,
at both panel widths (\cref{fig:sweep-length}). Since $N_{\mathrm{stages}}=s\,N_{\mathrm{sweeps}}$, the
sweep count depends little on the sweep length: in the median, the
shorter sweeps need $0.90$ to $0.94$ times as many sweeps as RR-OSBJ.
Shortening the sweep therefore reduces the stages to convergence in
about the same proportion.

The sweep count also depends little on the order of the stages and on
the panel width. Randomly permuting the stages of AG-OSBJ changes its
sweep count by at most $0.76$, and its sweep counts at $K=64$ and
$K=128$ differ by at most $1.9$ across the sizes studied below. This
supports choosing the panel width by the time per sweep
(\cref{sec:panel-width}).

We next examine whether this behavior persists across matrix sizes.
Since the reduction in sweep length is near its largest at powers of
two (\cref{fig:stages-relative}), we add three equally spaced sizes
between consecutive powers of two from $512$ to $8192$, giving $17$
sizes.

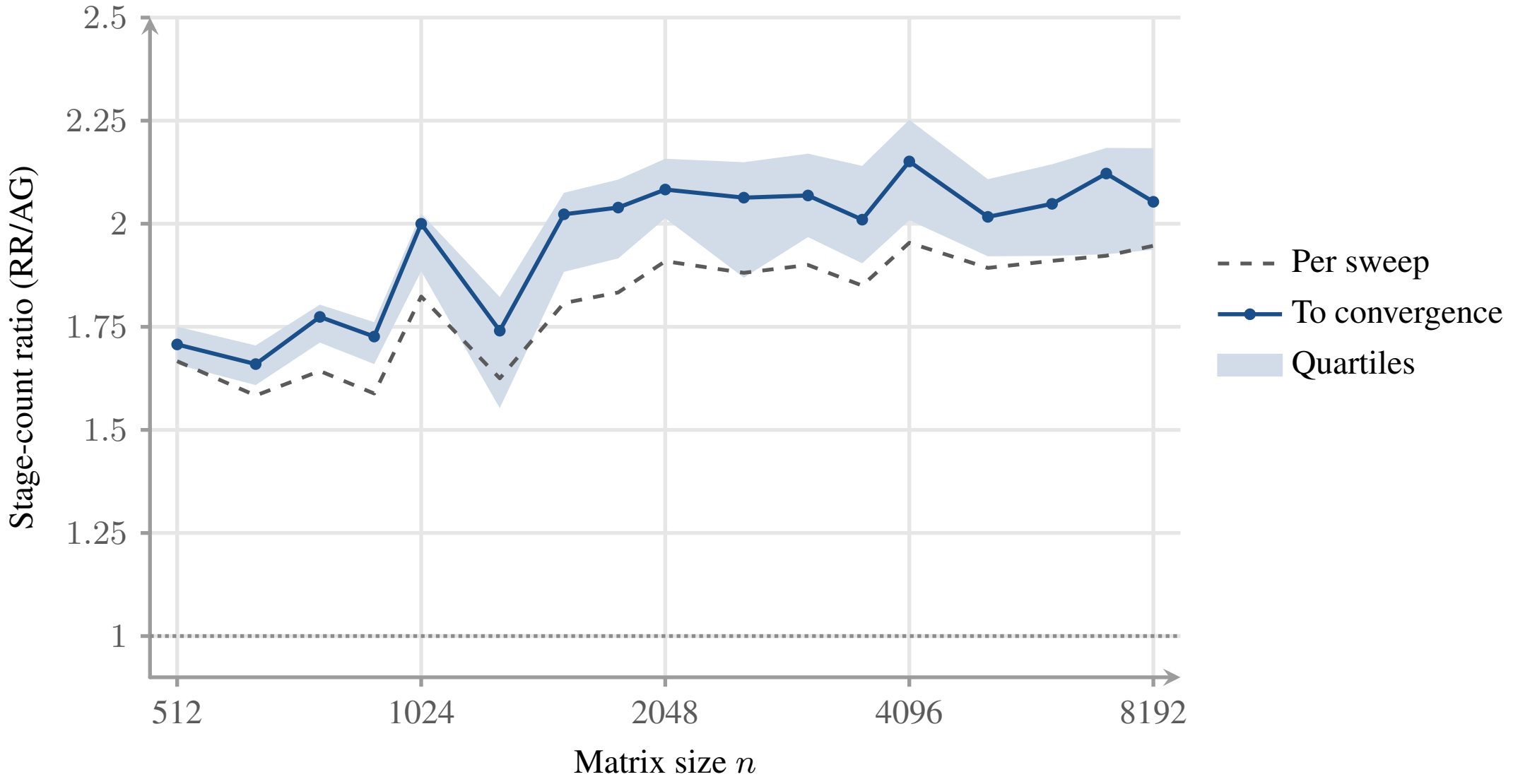
\begin{figure}[htbp]
\centering
\definecolor{ours}{HTML}{1B4F8A}
\definecolor{robin}{HTML}{C8553D}
\definecolor{bound}{HTML}{5A5A5A}
\definecolor{gesvdj}{HTML}{3A7D44}
\definecolor{gesvdp}{HTML}{C99A2E}
\begin{tikzpicture}
\begin{axis}[
    width=0.68\linewidth, height=0.47\linewidth,
    scale only axis=false,
    axis lines=left, axis line style={bound!60, line width=1bp, line cap=round},
    tick style={bound!60, line width=1bp}, tick align=outside,
    major tick length=2.5pt,
    ticklabel style={font=\footnotesize, text=bound},
    /pgf/number format/1000 sep={},
    label style={font=\small},
    every axis title/.style={font=\small, at={(0.5,1.04)}, anchor=south},
    grid=major, major grid style={bound!15, line width=1bp},
    every axis plot/.append style={line join=round},
    legend cell align=left,
    legend style={draw=none, fill=none,
        font=\scriptsize, row sep=-2pt, inner sep=1pt},
    xmode=log,
    log basis x=2,
    xmin=474,
    xmax=8847,
    xtick={512,1024,2048,4096,8192},
    xticklabels={512,1024,2048,4096,8192},
    xlabel={Matrix size $n$},
    ymin=0.9,
    ymax=2.50,
    ytick={1,1.25,1.5,1.75,2,2.25,2.5},
    ylabel={Stage-count ratio (RR/AG)},
    legend style={at={(1.03,0.55)}, anchor=west, font=\footnotesize, row sep=2pt},
]
\addplot[bound!70, densely dotted, line width=1bp, forget plot] coordinates {(474,1) (8847,1)};
\addplot[name path=upper, draw=none, forget plot] table[x=n,y=q3] {figures/stage_counts.dat};
\addplot[name path=lower, draw=none, forget plot] table[x=n,y=q1] {figures/stage_counts.dat};
\addplot[ours!20, forget plot] fill between[of=upper and lower];
\addplot[bound, dashed, line width=1bp] table[x=n,y=length] {figures/stage_counts.dat};
\addlegendentry{Per sweep}
\addplot[ours, line width=1.1bp, mark=*, mark size=1pt] table[x=n,y=median] {figures/stage_counts.dat};
\addlegendentry{To convergence}
\addlegendimage{area legend, draw=ours!20, fill=ours!20}
\addlegendentry{Quartiles}
\end{axis}
\end{tikzpicture}
\caption{Stages to convergence of RR-OSBJ divided by those of AG-OSBJ
against the matrix size $n$, at $K=64$: median and interquartile range
over the $39$ test matrices of each size. Dashed: the same ratio for the
stages per sweep.}
\label{fig:stages-ratio}
\end{figure}

\Cref{fig:stages-ratio} shows that the same holds at every size: the
reduction in total stages follows the reduction in stages per sweep.
At every size and panel width, the median ratio
$N^{\mathrm{RR}}_{\mathrm{stages}}/N^{\mathrm{AG}}_{\mathrm{stages}}$ is between $0.95$ and $1.12$ times
the ratio $s^{\mathrm{RR}}/s^{\mathrm{AG}}$. Run by run, both solvers
converge on every test matrix, and AG-OSBJ needs fewer stages than
RR-OSBJ in $1320$ of the $1326$ matched runs.

The exceptions come from two matrices. For \mfam{two plateaus}, shorter
sweeps consistently increase the sweep count: AG-OSBJ needs $1.11$ to
$1.70$ times as many sweeps as RR-OSBJ, depending on the size and
panel width, and at $n=640$ and $K=128$, where it uses $8$ stages per
sweep against $9$, it needs more stages on all three seeds. \mfam{One large}
with $\kappa=10^4$ and \mfam{one small} converge within four sweeps under
both solvers, so small differences in sweep count produce large
changes in their ratios; three runs on \mfam{one large} with $\kappa=10^4$
are the other exceptions.

The fewer stages also shorten the running time. We time the complete
solvers of \cref{sec:impl-overview}, which check the stopping
criterion after every sweep, as in \cref{sec:exp-setup}: on the $39$
test matrices at the $17$ sizes, with both solvers at the panel width
of \cref{sec:panel-width}, AG-OSBJ is $1.6$ to $2.0$ times faster than
RR-OSBJ in the median at each size. A stage of AG-OSBJ takes $4\%$ to $9\%$ longer than
one of RR-OSBJ, which offsets part of the reduction in stages.

\subsection{Comparison with dynamic ordering}
\label{sec:exp-dynamic}

We compare AG-OSBJ with the dynamic block ordering of
Be\v{c}ka, Ok\v{s}a, and Vajter\v{s}ic~\cite{BeckaOksaVajtersic2015}
and Kudo et al.~\cite{KudoEtAl2017}. This method retains a fixed
block partition but selects block pairs before every stage through
a greedy approximation to maximum-weight matching. We implement
it with the kernels of AG-OSBJ and call it DO-OSBJ (Dynamic-Ordering
OSBJ). Following
Be\v{c}ka et al., we assign blocks $i$ and $j$ the weight
$\|C_{ij}\|_F^2$, where $C_{ij}$ contains the cosines between
their columns. Pairs are considered in decreasing order of weight
and selected if both blocks are unmatched and at least one cosine
reaches the stopping tolerance. The remaining blocks are paired
in index order.

Following the published method, DO-OSBJ diagonalizes each
$H_S=X_S^\top X_S$ to convergence, repeating inner sweeps until a
complete sweep applies no rotation, and then sorts the columns of
$X_S$ by norm. From the second inner sweep on, off-diagonal entries
of the partially diagonalized $H_S$ of magnitude at most
$\epsmach\|H_S\|_F$ are treated as zero, since $H_S$ does not
determine them more accurately.

We compare the three solvers on the $39$ test matrices
of \cref{sec:exp-sweeps}, at the five powers of two from $512$ to $8192$
and both panel widths. As in that section, convergence is checked after every
stage. No solver uses the QR and LQ preprocessing of Kudo
et al.~\cite{KudoEtAl2017}, so all start from the same matrix.

\begin{figure}[htbp]
\centering
\definecolor{ours}{HTML}{1B4F8A}
\definecolor{robin}{HTML}{C8553D}
\definecolor{bound}{HTML}{5A5A5A}
\definecolor{gesvdj}{HTML}{3A7D44}
\definecolor{gesvdp}{HTML}{C99A2E}
\begin{tikzpicture}
\begin{axis}[
    width=0.68\linewidth, height=0.47\linewidth,
    scale only axis=false,
    axis lines=left, axis line style={bound!60, line width=1bp, line cap=round},
    tick style={bound!60, line width=1bp}, tick align=outside,
    major tick length=2.5pt,
    ticklabel style={font=\footnotesize, text=bound},
    /pgf/number format/1000 sep={},
    label style={font=\small},
    every axis title/.style={font=\small, at={(0.5,1.04)}, anchor=south},
    grid=major, major grid style={bound!15, line width=1bp},
    every axis plot/.append style={line join=round},
    legend cell align=left,
    legend style={draw=none, fill=none,
        font=\scriptsize, row sep=-2pt, inner sep=1pt},
    xmode=log,
    log basis x=2,
    xmin=457,
    xmax=9175,
    xtick={512,1024,2048,4096,8192},
    xticklabels={512,1024,2048,4096,8192},
    xlabel={Matrix size $n$},
    ymin=0.20,
    ymax=1.40,
    ytick={0.25,0.5,0.75,1,1.25,1.5},
    ylabel={Stage-count ratio (DO/AG)},
    legend style={at={(1.03,0.55)}, anchor=west, font=\footnotesize, row sep=2pt},
]
\addplot[bound!70, densely dotted, line width=1bp, forget plot] coordinates {(457,1) (9175,1)};
\addplot[ours, mark=*, line width=1.1bp, mark size=1.3pt, error bars/.cd, y dir=both, y explicit, error bar style={solid, line width=1bp}, error mark options={rotate=90, mark size=2.5pt, line width=1bp}] table[x=x64, y=median64, y error minus expr=\thisrow{median64}-\thisrow{q1_64}, y error plus expr=\thisrow{q3_64}-\thisrow{median64}] {figures/dynamic_ratio.dat};
\addlegendentry{$K=64$}
\addplot[robin, densely dashed, mark=square*, mark options={solid}, line width=1.1bp, mark size=1.3pt, error bars/.cd, y dir=both, y explicit, error bar style={solid, line width=1bp}, error mark options={rotate=90, mark size=2.5pt, line width=1bp}] table[x=x128, y=median128, y error minus expr=\thisrow{median128}-\thisrow{q1_128}, y error plus expr=\thisrow{q3_128}-\thisrow{median128}] {figures/dynamic_ratio.dat};
\addlegendentry{$K=128$}
\end{axis}
\end{tikzpicture}
\caption{Stages to convergence of DO-OSBJ divided by those of AG-OSBJ
against the matrix size $n$, at $K=64$ and $K=128$: median and
interquartile range over the $39$ test matrices at each size.}
\label{fig:dynamic-ratio}
\end{figure}
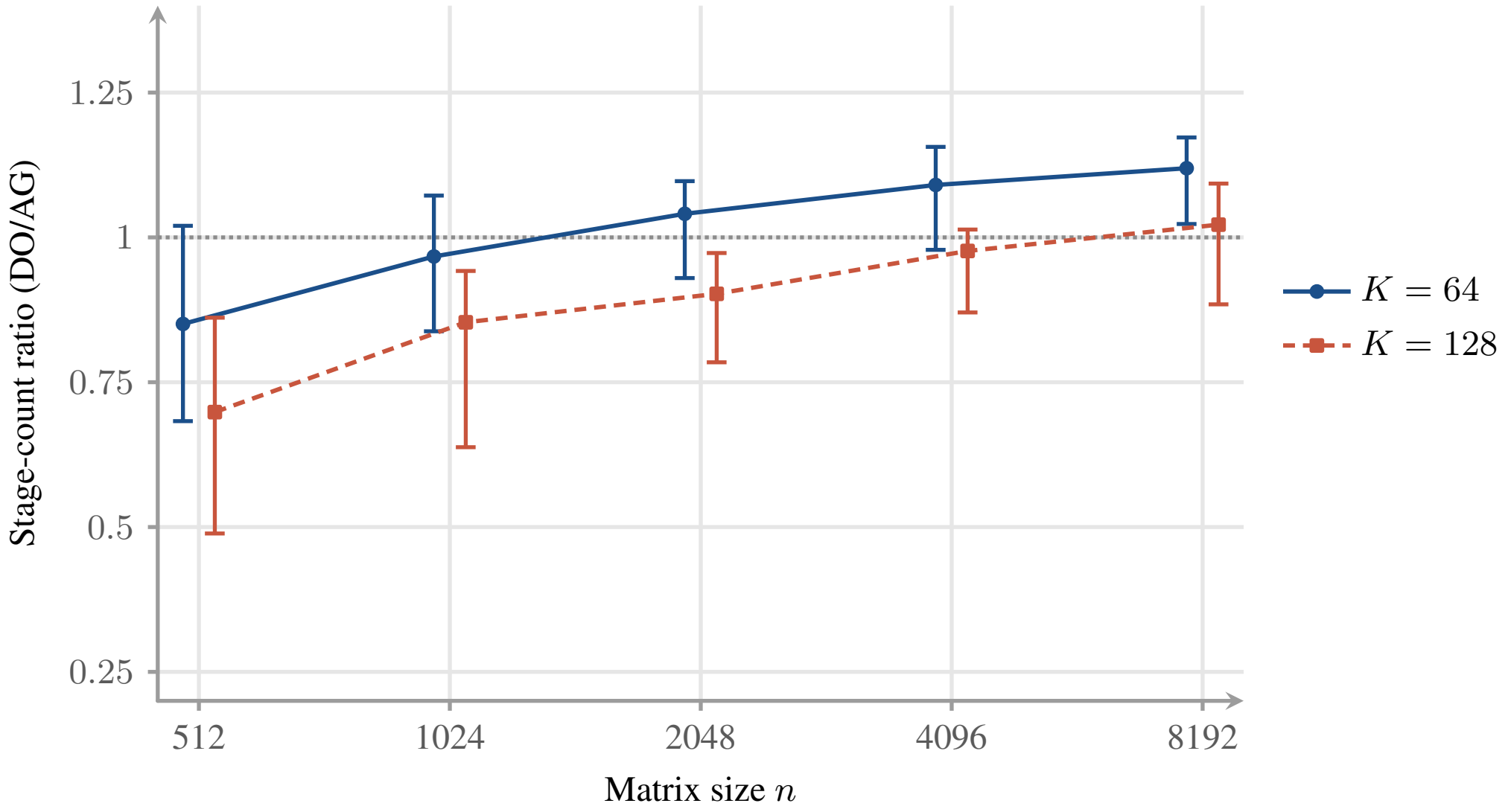

The advantage of DO-OSBJ in stages to convergence decreases with
matrix size (\cref{fig:dynamic-ratio}). At $K=64$, the median ratio
to AG-OSBJ rises from $0.85$ at $n=512$ to $1.12$ at $n=8192$;
at $K=128$, it rises from $0.70$ to $1.02$. Two matrices retain
an advantage for DO-OSBJ at every size, width, and seed:
\mfam{polynomial decay} and \mfam{two plateaus}.
\Cref{tab:dynamic-matrix} gives the results for individual
matrices at $n=2048$.

For a comparison with the same panel solve as AG-OSBJ, we also
run DO-OSBJ with one inner sweep per pivot set and no column
sorting. One sweep is not expected to suffice: dynamic ordering
pairs the blocks that are furthest from orthogonal, and
orthogonalizing them takes several inner sweeps. This
variant requires $1.29$ to $1.56$ times as many stages as AG-OSBJ
in the median, depending on the size and width.

Complete panel solves, however, increase the work within each stage.
Because pivot sets are processed concurrently, the next stage
must wait for the slowest panel. For stage $k$ with partition
$\cP_k$, let $r_{k,S}$ denote the number of inner sweeps on
pivot set $S$, including the final sweep that applies no rotation.
We measure the mean per-stage maximum
\[
    r_k^{\max}=\max_{S\in\cP_k}r_{k,S},
    \qquad
    \bar r_{\max}=\frac{1}{N_{\mathrm{stages}}}\sum_{k=1}^{N_{\mathrm{stages}}}r_k^{\max}.
\]

\begin{table}[t]
\caption{Stages to convergence of AG-OSBJ, RR-OSBJ, and DO-OSBJ
for each \texttt{randsvd} matrix at $n=2048$ (first seed).
Bold: fewest stages at each panel width.
Parentheses: mean per-stage maximum inner-sweep count
$\bar r_{\max}$ for DO-OSBJ.}
\label{tab:dynamic-matrix}
\centering
\footnotesize
\setlength{\tabcolsep}{4pt}
\begin{tabular}{@{}llrrrrrr@{}}
\toprule
& & \multicolumn{3}{c}{$K=64$} & \multicolumn{3}{c}{$K=128$} \\
\cmidrule(lr){3-5} \cmidrule(l){6-8}
Matrix & $\kappa$ & AG & RR & DO & AG & RR & DO \\
\midrule
\mfam{polynomial decay} & $n^{2}$ & 784 & 1637 & \textbf{609} (6.3) & 394 & 737 & \textbf{267} (6.6) \\
\mfam{two plateaus} & $10^{4}$ & 747 & 1130 & \textbf{430} (5.0) & 379 & 564 & \textbf{225} (6.3) \\
\mfam{arithmetic} & $10^{4}$ & \textbf{451} & 950 & 462 (5.7) & 229 & 449 & \textbf{216} (5.8) \\
\mfam{arithmetic} & $10^{10}$ & \textbf{439} & 944 & 495 (5.5) & 237 & 464 & \textbf{226} (5.7) \\
\mfam{log-random} & $10^{4}$ & \textbf{754} & 1603 & 785 (6.5) & 404 & 806 & \textbf{345} (6.8) \\
\mfam{log-random} & $10^{10}$ & 1466 & 3233 & \textbf{1365} (6.8) & 766 & 1598 & \textbf{606} (7.2) \\
\mfam{geometric} & $10^{4}$ & \textbf{729} & 1568 & 781 (6.5) & 381 & 782 & \textbf{347} (6.8) \\
\mfam{geometric} & $10^{10}$ & 1446 & 3274 & \textbf{1380} (6.8) & 754 & 1586 & \textbf{590} (7.3) \\
\mfam{multi-cluster} & $n$ & 474 & 910 & \textbf{465} (5.7) & 251 & 457 & \textbf{213} (5.8) \\
\mfam{one large} & $10^{4}$ & 29 & 57 & \textbf{12} (5.8) & 17 & 38 & \textbf{8} (7.6) \\
\mfam{one large} & $10^{10}$ & \textbf{345} & 729 & 384 (5.5) & 183 & 363 & \textbf{182} (5.7) \\
\mfam{one small} & $10^{4}$ & \textbf{53} & 121 & 105 (2.0) & \textbf{29} & 60 & 50 (2.0) \\
\mfam{one small} & $10^{10}$ & \textbf{53} & 121 & 110 (2.0) & \textbf{29} & 60 & 48 (2.0) \\
\midrule
\multicolumn{2}{@{}l}{Median ratio to AG} & 1 & 2.13 & 1.02 & 1 & 2.00 & 0.85 \\
\bottomrule
\end{tabular}

\end{table}

Across the tested sizes and widths, the median $\bar r_{\max}$
ranges from $5.4$ to $6.7$, against one inner sweep per stage
for AG-OSBJ. Summed to convergence, DO-OSBJ therefore waits
for $3.8$ to $6.4$ times as many inner sweeps in the median.

Dynamic ordering also adds work before every stage, whatever the
panel solve. The weights depend on the current iterate, so they
require a separate pass over the whole matrix. The exact weights
used here form all inner products between columns, about $n/K$
times the arithmetic of \kgram{} in one stage. Greedy matching then selects
the $n/K$ pivot sets one after another, a chain of dependent steps
whose cost is set by latency. AG-OSBJ uses a precomputed schedule
and needs neither step.

\subsection{Accuracy and running time against cuSOLVER}
\label{sec:exp-cusolver}

We finally assess the numerical accuracy of AG-OSBJ and compare
its relative accuracy and running time with cuSOLVER's
double-precision \texttt{gesvdj}, \texttt{gesvd}, and
\texttt{gesvdp} routines~\cite{cuSOLVER}.
These routines are based on Jacobi iterations,
bidiagonalization followed by QR iteration, and polar
decomposition followed by a symmetric eigensolver, respectively.
All methods compute full $U$ and $V$.
The experiments use the first seed of every matrix and
$n\in\{512,1024,2048,4096,8192\}$.
AG-OSBJ runs as described in \cref{sec:implementation},
with a limit of $100$ outer sweeps, which no run reached.

\paragraph{Numerical accuracy}
Let $\sigma^*=(\sigma^*_i)$ and $\hat\sigma=(\hat\sigma_i)$
denote the reference and computed singular values,
respectively, in decreasing order.
For the \texttt{randsvd} matrices, the reference values are the
prescribed $\sigma_i$.
The largest errors of AG-OSBJ on these matrices, over all sizes
and both panel widths, are as follows:
\begin{itemize}[leftmargin=*]
\item backward error
  $\eta=\|A-\hat U\hat\Sigma\hat V^\top\|_F/(n\|A\|_F)$,
  at most $1.3\cdot10^{-17}$;
\item orthogonality errors
  $\rho_U=\|I-\hat U^\top\hat U\|_F/n$ and
  $\rho_V=\|I-\hat V^\top\hat V\|_F/n$, at most
  $1.1\cdot10^{-14}$ (about $49\epsmach$) and
  $5.9\cdot10^{-16}$, respectively;
\item singular-value error
  $\rho_S=\|\hat\sigma-\sigma^*\|_2/(n\|\sigma^*\|_2)$,
  at most $1.8\cdot10^{-18}$.
\end{itemize}
The orthogonality of $\hat U$ is the least accurate. As pairs of
columns whose cosine is at most $\sqrt m\,\epsmach$ in absolute value
are not rotated (\cref{sec:impl-overview}), the cosines between the
final columns of $X$, and so between those of $\hat U$, remain at
about that level.
A similar limitation is observed in \texttt{gesvdj}, for
which the largest $\rho_U$ on the same matrices is
$9.3\cdot10^{-14}$, about eight times that of AG-OSBJ.

\paragraph{Relative accuracy}
For the column-graded matrices, the reference singular values are
computed by LAPACK's \texttt{DGEJSV}, a preconditioned Jacobi SVD
routine designed to achieve high relative accuracy on such
matrices~\cite{DrmacVeselic2008a,DrmacVeselic2008b}.\footnote{
We use \texttt{JOBA='G'}, which enables full pivoting and
condition estimation.}
Since a normwise error can conceal inaccuracies in individual
small singular values, we measure the largest relative
singular-value error
\[
\rho_{S,\max}=
\max_i\frac{|\hat\sigma_i-\sigma^*_i|}{\sigma^*_i}.
\]

\begin{table}[htbp]
\caption{Largest relative singular-value error $\rho_{S,\max}$
on the column-graded matrices at each $\kappa_D$, over the
powers of two $n=512,\dots,8192$, of AG-OSBJ (either panel
width) and three cuSOLVER routines.
Bold: the smallest in each column.}
\label{tab:accuracy}
\centering
\footnotesize
\setlength{\tabcolsep}{4pt}
\begin{tabular}{@{}lcc@{}}
\toprule
Solver & $\kappa_D=10^{4}$ & $\kappa_D=10^{10}$ \\
\midrule
AG-OSBJ & {\boldmath$1.5\cdot10^{-13}$} & {\boldmath$8.0\cdot10^{-14}$} \\
\texttt{gesvdj} & $5.9\cdot10^{-12}$ & $2.9\cdot10^{-8}$ \\
\addlinespace
\texttt{gesvd} & $2.4\cdot10^{-13}$ & $2.0\cdot10^{-7}$ \\
\texttt{gesvdp} & $5.9\cdot10^{-13}$ & $2.4\cdot10^{-6}$ \\
\bottomrule
\end{tabular}

\end{table}

As shown in \cref{tab:accuracy}, AG-OSBJ agrees with
\texttt{DGEJSV} to $O(n\epsmach)$ at every size,
independently of $\kappa_D$.
In contrast, \texttt{gesvdj} loses relative accuracy at
$\kappa_D=10^{10}$, where its largest $\rho_{S,\max}$
exceeds that of AG-OSBJ by more than five orders of magnitude.
Its normwise stopping criterion does not enforce the
columnwise criterion \eqref{eq:stopping} required for
relative accuracy here~\cite{DemmelVeselic1992}.
The general-purpose routines \texttt{gesvd} and
\texttt{gesvdp} likewise lose relative accuracy on strongly
column-graded matrices, despite their normwise backward
stability.

\paragraph{Running time}
The timings use the \texttt{randsvd} matrices at the powers of two
above and, to assess performance at sizes that are not multiples of
either panel width, at $n\in\{724,1448,2896,5793\}$.
AG-OSBJ uses the panel width selected in
\cref{sec:panel-width}, and its timings include schedule
construction and CUDA graph capture.

\begin{table}[htbp]
\caption{Speedup of AG-OSBJ over cuSOLVER.
For each \texttt{randsvd} matrix and size $n$: time of
\texttt{gesvdp}, the fastest cuSOLVER routine in all but
two tests, divided by that of AG-OSBJ.
Bottom rows: smallest of the same ratio over the $13$
matrices for \texttt{gesvdj} and \texttt{gesvd}.
Values above $1$, in bold, mean AG-OSBJ is faster.}
\label{tab:timing}
\centering
\footnotesize
\setlength{\tabcolsep}{4pt}
\renewcommand{\arraystretch}{1.1}
\begin{tabular}{@{}llrrrrrrrrr@{}}
\toprule
& & \multicolumn{9}{c}{$n$} \\
\cmidrule(l){3-11}
Spectrum & $\kappa$ & 512 & 724 & 1024 & 1448 & 2048 & 2896 & 4096 & 5793 & 8192 \\
\midrule
\mfam{polynomial decay} & $n^{2}$ & \textbf{1.68} & \textbf{1.52} & \textbf{1.69} & \textbf{1.29} & \textbf{1.15} & 0.76 & 0.83 & 0.58 & 0.53 \\
\addlinespace[2pt]
\mfam{two plateaus} & $10^{4}$ & \textbf{1.40} & \textbf{1.28} & \textbf{1.56} & \textbf{1.25} & \textbf{1.16} & 0.85 & 0.95 & 0.73 & 0.68 \\
\addlinespace[2pt]
\multirow{2}{*}{\mfam{arithmetic}} & $10^{4}$ & \textbf{2.34} & \textbf{2.51} & \textbf{2.91} & \textbf{2.28} & \textbf{1.97} & \textbf{1.36} & \textbf{1.43} & \textbf{1.10} & 0.93 \\
 & $10^{10}$ & \textbf{2.77} & \textbf{2.48} & \textbf{2.83} & \textbf{2.41} & \textbf{2.05} & \textbf{1.33} & \textbf{1.40} & \textbf{1.08} & 0.91 \\
\addlinespace[2pt]
\multirow{2}{*}{\mfam{log-random}} & $10^{4}$ & \textbf{1.54} & \textbf{1.38} & \textbf{1.62} & \textbf{1.29} & \textbf{1.20} & 0.80 & 0.89 & 0.66 & 0.57 \\
 & $10^{10}$ & 0.88 & 0.80 & 0.89 & 0.74 & 0.65 & 0.44 & 0.47 & 0.36 & 0.31 \\
\addlinespace[2pt]
\multirow{2}{*}{\mfam{geometric}} & $10^{4}$ & \textbf{1.60} & \textbf{1.44} & \textbf{1.68} & \textbf{1.34} & \textbf{1.20} & 0.83 & 0.89 & 0.68 & 0.59 \\
 & $10^{10}$ & 0.90 & 0.82 & 0.89 & 0.73 & 0.65 & 0.44 & 0.47 & 0.36 & 0.31 \\
\addlinespace[2pt]
\mfam{multi-cluster} & $n$ & \textbf{2.03} & \textbf{2.11} & \textbf{2.15} & \textbf{1.77} & \textbf{1.56} & \textbf{1.20} & \textbf{1.43} & 0.97 & 0.84 \\
\addlinespace[2pt]
\multirow{2}{*}{\mfam{one large}} & $10^{4}$ & \textbf{6.78} & \textbf{12.07} & \textbf{14.49} & \textbf{20.74} & \textbf{20.71} & \textbf{14.05} & \textbf{10.00} & \textbf{14.90} & \textbf{14.93} \\
 & $10^{10}$ & \textbf{2.96} & \textbf{3.04} & \textbf{3.50} & \textbf{2.71} & \textbf{2.40} & \textbf{1.74} & \textbf{1.81} & \textbf{1.40} & \textbf{1.26} \\
\addlinespace[2pt]
\multirow{2}{*}{\mfam{one small}} & $10^{4}$ & \textbf{11.77} & \textbf{12.15} & \textbf{14.07} & \textbf{11.98} & \textbf{11.88} & \textbf{8.08} & \textbf{10.14} & \textbf{7.82} & \textbf{7.76} \\
 & $10^{10}$ & \textbf{12.35} & \textbf{12.96} & \textbf{14.53} & \textbf{13.03} & \textbf{12.34} & \textbf{8.48} & \textbf{10.60} & \textbf{8.19} & \textbf{8.04} \\
\midrule
\multicolumn{2}{@{}l}{\emph{Smallest} over \texttt{gesvdj}} & \textbf{3.88} & \textbf{3.74} & \textbf{4.86} & \textbf{5.23} & \textbf{5.34} & \textbf{4.70} & \textbf{5.29} & \textbf{5.38} & \textbf{4.90} \\
\multicolumn{2}{@{}l}{\emph{Smallest} over \texttt{gesvd}} & \textbf{2.06} & \textbf{1.67} & \textbf{1.88} & \textbf{1.62} & \textbf{1.46} & \textbf{1.01} & \textbf{1.11} & 0.73 & 0.79 \\
\bottomrule
\end{tabular}

\end{table}

AG-OSBJ is faster than the Jacobi routine \texttt{gesvdj} in
every test, by a factor of at least $3.7$ and of up to $321$ on
\mfam{one large} with $\kappa=10^4$.
Excluding the \mfam{one large} and \mfam{one small}
families, the median speedup over the remaining nine matrices
is $6.0$.
AG-OSBJ is also faster than \texttt{gesvd} in all but two
tests (\mfam{log-random} with $\kappa=10^{10}$ at
$n=5793$ and $8192$), with a median speedup of $4.3$
over the same nine matrices.

The comparison with \texttt{gesvdp}, the fastest cuSOLVER
routine in all but two tests, is more dependent on the
matrix size and spectrum.
AG-OSBJ is faster on most matrices up to $n=2048$,
on about half at intermediate sizes
($2896\leq n\leq5793$), and on the \mfam{one large}
and \mfam{one small} families at every size.
Conversely, \texttt{gesvdp} is faster on the two spectra
requiring the most sweeps, \mfam{geometric} and
\mfam{log-random} with $\kappa=10^{10}$, and at $n=8192$
on all matrices outside the \mfam{one large} and
\mfam{one small} families.

\renewcommand\bottomfraction{.1}
\Needspace{9\baselineskip}
\section{Conclusion}
\label{sec:conclusion}
This paper has shown that variable partitioning can substantially
reduce the sequential depth of the one-sided block Jacobi SVD
without sacrificing the computational advantages of blocking.
Using affine planes and transversal designs, we construct static
schedules for arbitrary $n$ and $K$ that nearly halve the number
of stages per sweep compared with block round-robin and attain
the combinatorial lower bound for several parameter ranges.

The experiments show that this reduction in sequential depth
translates into shorter running times. Shorter sweeps
generally do not require additional sweeps to converge, with
AG-OSBJ using fewer total stages than RR-OSBJ in $1320$ of
$1326$ matched runs. On an NVIDIA H200 NVL, the resulting
median speedup over RR-OSBJ ranges from $1.6$ to $2.0$ across the $17$
tested sizes between $512$ and $8192$.
Static scheduling also provides a favorable alternative to
dynamic ordering: although the latter converges in fewer
stages at small sizes, this advantage disappears as $n$ grows,
while its higher inner iteration counts and per-stage
scheduling costs remain.

AG-OSBJ is also competitive with established GPU SVD
routines. It outperforms
cuSOLVER's \texttt{gesvdj} in every test and \texttt{gesvd}
in all but two. Against \texttt{gesvdp}, it is faster on most
matrices up to $n=2048$, although \texttt{gesvdp} is generally
faster at $n=8192$.
These performance gains are achieved while maintaining
numerical errors close to roundoff level and high relative
accuracy in the singular values of strongly column-graded
matrices, which all three cuSOLVER routines lose.

AG-OSBJ converged on every tested matrix. In exact
arithmetic, it belongs to the class of quasi-cyclic Jacobi
methods, so its asymptotic convergence is quadratic when
the singular values are simple.
Establishing its global convergence under variable
partitioning and approximate panel solves remains an
important direction for future work.
\newcommand{\paperbackmatter}{%

    \pdfbookmark[1]{Acknowledgments}{preprint-acknowledgments}

  \section*{Acknowledgments}
  The authors used AI tools throughout this work, including for writing
  and testing code, running and analyzing experiments, preparing figures,
  checking derivations, and writing and revising the manuscript. The
  research ideas and the design of the experiments are the authors' own,
  as is the interpretation of the results. The authors take full
  responsibility for all content.

    \pdfbookmark[1]{References}{preprint-references}

  \bibliographystyle{siamplain}
  \bibliography{references}}

\begin{thebibliography}{10}

\bibitem{AbdelfattahFasi2026}
{\sc A.~Abdelfattah and M.~Fasi}, {\em An efficient batch solver for the
  singular value decomposition on {GPUs}}.
\newblock arXiv:2601.17979, 2026, \url{https://arxiv.org/abs/2601.17979}.

\bibitem{AbelColbournDinitz2007}
{\sc R.~J.~R. Abel, C.~J. Colbourn, and J.~H. Dinitz}, {\em Mutually orthogonal
  {L}atin squares ({MOLS})}, in Handbook of Combinatorial Designs, C.~J.
  Colbourn and J.~H. Dinitz, eds., Chapman \& Hall/CRC, Boca Raton, FL,
  2nd~ed., 2007, ch.~III.3, pp.~160--193.

\bibitem{BakerHarmanPintz2001}
{\sc R.~C. Baker, G.~Harman, and J.~Pintz}, {\em The difference between
  consecutive primes, {II}}, Proc. London Math. Soc. (3), 83 (2001),
  pp.~532--562, \url{https://doi.org/10.1112/plms/83.3.532}.

\bibitem{BeckaOksa2016}
{\sc M.~Be{\v c}ka and G.~Ok{\v s}a}, {\em New approach to local computations
  in the parallel one-sided {Jacobi} {SVD} algorithm}, in Parallel Processing
  and Applied Mathematics, R.~Wyrzykowski et~al., eds., vol.~9573 of Lecture
  Notes in Comput. Sci., Cham, 2016, Springer, pp.~605--617,
  \url{https://doi.org/10.1007/978-3-319-32149-3_56}.

\bibitem{BeckaOksaVajtersic2002}
{\sc M.~Be{\v c}ka, G.~Ok{\v s}a, and M.~Vajter{\v s}ic}, {\em Dynamic ordering
  for a parallel block-{Jacobi} {SVD} algorithm}, Parallel Comput., 28 (2002),
  pp.~243--262, \url{https://doi.org/10.1016/S0167-8191(01)00138-7}.

\bibitem{BeckaOksaVajtersic2015}
{\sc M.~Be{\v c}ka, G.~Ok{\v s}a, and M.~Vajter{\v s}ic}, {\em New dynamic
  orderings for the parallel one-sided block-{Jacobi} {SVD} algorithm},
  Parallel Process. Lett., 25 (2015), 1550003,
  \url{https://doi.org/10.1142/S0129626415500036}.

\bibitem{BeckaVajtersic1999}
{\sc M.~Be{\v c}ka and M.~Vajter{\v s}ic}, {\em Block-{Jacobi} {SVD} algorithms
  for distributed memory systems {I}: Hypercubes and rings}, Parallel
  Algorithms Appl., 13 (1999), pp.~265--287,
  \url{https://doi.org/10.1080/10637199808947370}.

\bibitem{Bierbrauer1991}
{\sc J.~Bierbrauer}, {\em Weighted arcs, the finite {Radon} transform and a
  {Ramsey} problem}, Graphs Combin., 7 (1991), pp.~113--118,
  \url{https://doi.org/10.1007/BF01788137}.

\bibitem{BierbrauerGyarfas1987}
{\sc J.~Bierbrauer and A.~Gy\'arf\'as}, {\em On $(n,k)$-colorings of complete
  graphs}, Congr. Numer., 58 (1987), pp.~123--139.

\bibitem{BoukaramEtAl2018}
{\sc W.~H. Boukaram, G.~Turkiyyah, H.~Ltaief, and D.~E. Keyes}, {\em Batched
  {QR} and {SVD} algorithms on {GPUs} with applications in hierarchical matrix
  compression}, Parallel Comput., 74 (2018), pp.~19--33,
  \url{https://doi.org/10.1016/j.parco.2017.09.001}.

\bibitem{BrentLuk1985}
{\sc R.~P. Brent and F.~T. Luk}, {\em The solution of singular-value and
  symmetric eigenvalue problems on multiprocessor arrays}, SIAM J. Sci.
  Statist. Comput., 6 (1985), pp.~69--84,
  \url{https://doi.org/10.1137/0906007}.

\bibitem{ColbournDinitz2007}
{\sc C.~J. Colbourn and J.~H. Dinitz}, eds., {\em Handbook of Combinatorial
  Designs}, Discrete Mathematics and its Applications, Chapman \& Hall/CRC,
  Boca Raton, FL, 2nd~ed., 2007, \url{https://doi.org/10.1201/9781420010541}.

\bibitem{CuencaGimenez2000}
{\sc J.~Cuenca and D.~Gim{\'e}nez}, {\em Implementation of parallel one-sided
  block {Jacobi} methods for the symmetric eigenvalue problem}, in Parallel
  Computing: Fundamentals and Applications (ParCo99), Imperial College Press,
  2000, pp.~291--298.

\bibitem{deRijk1989}
{\sc P.~P.~M. de~Rijk}, {\em A one-sided {Jacobi} algorithm for computing the
  singular value decomposition on a vector computer}, SIAM J. Sci. Statist.
  Comput., 10 (1989), pp.~359--371, \url{https://doi.org/10.1137/0910023}.

\bibitem{DemmelVeselic1992}
{\sc J.~Demmel and K.~Veseli{\'c}}, {\em Jacobi's method is more accurate than
  {QR}}, SIAM J. Matrix Anal. Appl., 13 (1992), pp.~1204--1245,
  \url{https://doi.org/10.1137/0613074}.

\bibitem{Drmac1997}
{\sc Z.~Drma{\v c}}, {\em Implementation of {Jacobi} rotations for accurate
  singular value computation in floating point arithmetic}, SIAM J. Sci.
  Comput., 18 (1997), pp.~1200--1222,
  \url{https://doi.org/10.1137/S1064827594265095}.

\bibitem{Drmac2009}
{\sc Z.~Drma{\v c}}, {\em A global convergence proof for cyclic {Jacobi}
  methods with block rotations}, SIAM J. Matrix Anal. Appl., 31 (2009),
  pp.~1329--1350, \url{https://doi.org/10.1137/090748548}.

\bibitem{DrmacVeselic2008a}
{\sc Z.~Drma{\v c} and K.~Veseli{\'c}}, {\em New fast and accurate {Jacobi}
  {SVD} algorithm. {I}}, SIAM J. Matrix Anal. Appl., 29 (2008), pp.~1322--1342,
  \url{https://doi.org/10.1137/050639193}.

\bibitem{DrmacVeselic2008b}
{\sc Z.~Drma{\v c} and K.~Veseli{\'c}}, {\em New fast and accurate {Jacobi}
  {SVD} algorithm. {II}}, SIAM J. Matrix Anal. Appl., 29 (2008),
  pp.~1343--1362, \url{https://doi.org/10.1137/05063920X}.

\bibitem{Eberlein1987}
{\sc P.~J. Eberlein}, {\em On one-sided {Jacobi} methods for parallel
  computation}, SIAM J. Algebraic Discrete Methods, 8 (1987), pp.~790--796,
  \url{https://doi.org/10.1137/0608064}.

\bibitem{EberleinPark1990}
{\sc P.~J. Eberlein and H.~Park}, {\em Efficient implementation of {Jacobi}
  algorithms and {Jacobi} sets on distributed memory architectures}, J.
  Parallel Distrib. Comput., 8 (1990), pp.~358--366,
  \url{https://doi.org/10.1016/0743-7315(90)90134-B}.

\bibitem{ForsytheHenrici1960}
{\sc G.~E. Forsythe and P.~Henrici}, {\em The cyclic {Jacobi} method for
  computing the principal values of a complex matrix}, Trans. Amer. Math. Soc.,
  94 (1960), pp.~1--23,
  \url{https://doi.org/10.1090/S0002-9947-1960-0109825-2}.

\bibitem{Furedi1981}
{\sc Z.~F\"uredi}, {\em Maximum degree and fractional matchings in uniform
  hypergraphs}, Combinatorica, 1 (1981), pp.~155--162,
  \url{https://doi.org/10.1007/BF02579271}.

\bibitem{Furedi1989}
{\sc Z.~F\"uredi}, {\em Covering the complete graph by partitions}, Discrete
  Math., 75 (1989), pp.~217--226,
  \url{https://doi.org/10.1016/0012-365X(89)90088-5}.

\bibitem{FurediGyarfas1991}
{\sc Z.~F\"uredi and A.~Gy\'arf\'as}, {\em Covering $t$-element sets by
  partitions}, European J. Combin., 12 (1991), pp.~483--489,
  \url{https://doi.org/10.1016/S0195-6698(13)80099-6}.

\bibitem{GolubVanLoan2013}
{\sc G.~H. Golub and C.~F. {Van Loan}}, {\em Matrix Computations}, Johns
  Hopkins University Press, Baltimore, MD, 4th~ed., 2013.

\bibitem{Gyarfas1977}
{\sc A.~Gy\'arf\'as}, {\em Partition Covers and Blocking Sets in Hypergraphs},
  PhD thesis, MTA SZTAKI, Budapest, 1977.
\newblock Tanulm\'anyok 71. In Hungarian.

\bibitem{Gyarfas2010}
{\sc A.~Gy\'arf\'as}, {\em Large monochromatic components in edge colorings of
  graphs: A survey}, in Ramsey Theory: Yesterday, Today, and Tomorrow,
  A.~Soifer, ed., vol.~285 of Progr. Math., Birkh\"auser, Boston, 2011,
  pp.~77--96, \url{https://doi.org/10.1007/978-0-8176-8092-3_5}.

\bibitem{Haemers1999}
{\sc W.~H. Haemers}, {\em Minimum resolvable coverings with small parallel
  classes}, Discrete Math., 197--198 (1999), pp.~393--396,
  \url{https://doi.org/10.1016/S0012-365X(99)90089-4}.

\bibitem{Hari2007}
{\sc V.~Hari}, {\em Quadratic convergence of a special quasi-cyclic {Jacobi}
  method}, Ann. Univ. Ferrara, 53 (2007), pp.~255--269,
  \url{https://doi.org/10.1007/s11565-007-0027-8}.

\bibitem{Hari2025}
{\sc V.~Hari}, {\em Global and quadratic convergence of the block {Jacobi}
  method for {H}ermitian matrices under the de {R}ijk pivot strategy},
  Electron. Trans. Numer. Anal., 63 (2025), pp.~83--128,
  \url{https://doi.org/10.1553/etna_vol63s83}.

\bibitem{HariBegovic2017}
{\sc V.~Hari and E.~{Begovi{\'c} Kova{\v c}}}, {\em Convergence of the cyclic
  and quasi-cyclic block {Jacobi} methods}, Electron. Trans. Numer. Anal., 46
  (2017), pp.~107--147,
  \url{https://etna.ricam.oeaw.ac.at/vol.46.2017/pp107-147.dir/pp107-147.pdf}.

\bibitem{Hestenes1958}
{\sc M.~R. Hestenes}, {\em Inversion of matrices by biorthogonalization and
  related results}, J. Soc. Indust. Appl. Math., 6 (1958), pp.~51--90,
  \url{https://doi.org/10.1137/0106005}.

\bibitem{Higham2002}
{\sc N.~J. Higham}, {\em Accuracy and Stability of Numerical Algorithms}, SIAM,
  Philadelphia, second~ed., 2002,
  \url{https://doi.org/10.1137/1.9780898718027}.

\bibitem{Hirschfeld1998}
{\sc J.~W.~P. Hirschfeld}, {\em Projective Geometries over Finite Fields},
  Oxford Math. Monogr., Clarendon Press, Oxford University Press, Oxford,
  2nd~ed., 1998, \url{https://doi.org/10.1093/oso/9780198502951.001.0001}.

\bibitem{HughesPiper1985}
{\sc D.~R. Hughes and F.~C. Piper}, {\em Design Theory}, Cambridge University
  Press, Cambridge, 1985.

\bibitem{Kaiser1972}
{\sc H.~F. Kaiser}, {\em The {JK} method: A procedure for finding the
  eigenvectors and eigenvalues of a real symmetric matrix}, Comput. J., 15
  (1972), pp.~271--273, \url{https://doi.org/10.1093/comjnl/15.3.271}.

\bibitem{KudoEtAl2017}
{\sc S.~Kudo, Y.~Yamamoto, M.~Be{\v c}ka, and M.~Vajter{\v s}ic}, {\em
  Performance analysis and optimization of the parallel one-sided block
  {Jacobi} {SVD} algorithm with dynamic ordering and variable blocking},
  Concurrency Computat.: Pract. Exper., 29 (2017), e4059,
  \url{https://doi.org/10.1002/cpe.4059}.

\bibitem{LukPark1989}
{\sc F.~T. Luk and H.~Park}, {\em On parallel {Jacobi} orderings}, SIAM J. Sci.
  Statist. Comput., 10 (1989), pp.~18--26,
  \url{https://doi.org/10.1137/0910002}.

\bibitem{LukPark1989proof}
{\sc F.~T. Luk and H.~Park}, {\em A proof of convergence for two parallel
  {Jacobi} {SVD} algorithms}, IEEE Trans. Comput., 38 (1989), pp.~806--811,
  \url{https://doi.org/10.1109/12.24289}.

\bibitem{Nagura1952}
{\sc J.~Nagura}, {\em On the interval containing at least one prime number},
  Proc. Japan Acad., 28 (1952), pp.~177--181,
  \url{https://doi.org/10.3792/pja/1195570997}.

\bibitem{Novakovic2015}
{\sc V.~Novakovi{\'c}}, {\em A hierarchically blocked {Jacobi} {SVD} algorithm
  for single and multiple graphics processing units}, SIAM J. Sci. Comput., 37
  (2015), pp.~C1--C30, \url{https://doi.org/10.1137/140952429},
  \url{https://arxiv.org/abs/1401.2720}.

\bibitem{Novakovic2023}
{\sc V.~Novakovi{\'c}}, {\em Vectorization of a thread-parallel {Jacobi}
  singular value decomposition method}, SIAM J. Sci. Comput., 45 (2023),
  pp.~C73--C100, \url{https://doi.org/10.1137/22M1478847},
  \url{https://arxiv.org/abs/2202.08361}.

\bibitem{CUTLASS}
{\sc {NVIDIA}}, {\em {CUTLASS}: {CUDA} templates for linear algebra
  subroutines}.
\newblock \url{https://github.com/NVIDIA/cutlass}, 2026.

\bibitem{cuSOLVER}
{\sc {NVIDIA Corporation}}, {\em {cuSOLVER} library documentation, release
  12.9}.
\newblock
  \url{https://docs.nvidia.com/cuda/archive/12.9.0/cusolver/index.html}.
\newblock CUDA Toolkit 12.9; accessed 2026-09-30.

\bibitem{OksaBecka2025}
{\sc G.~Ok{\v s}a and M.~Be{\v c}ka}, {\em Block-{Jacobi} {SVD} algorithms: A
  review}, Comput. Inform., 44 (2025), pp.~1254--1272,
  \url{https://doi.org/10.31577/cai_2025_5_1254}.

\bibitem{OksaVajtersic2006}
{\sc G.~Ok{\v s}a and M.~Vajter{\v s}ic}, {\em New clique-based parallel
  orderings for the block-{Jacobi} {EVD/SVD} algorithm}, Tech. Report 2006-05,
  Department of Computer Sciences, University of Salzburg, Dec. 2006,
  \url{https://www.cosy.sbg.ac.at/research/tr/2006-05_Oksa_Vajtersic.pdf}.

\bibitem{OksaYamamotoBeckaVajtersic2018}
{\sc G.~Ok{\v s}a, Y.~Yamamoto, M.~Be{\v c}ka, and M.~Vajter{\v s}ic}, {\em
  Asymptotic quadratic convergence of the parallel block-{Jacobi} {EVD}
  algorithm with dynamic ordering for {H}ermitian matrices}, BIT, 58 (2018),
  pp.~1099--1123, \url{https://doi.org/10.1007/s10543-018-0711-3}.

\bibitem{OksaYamamotoBeckaVajtersic2019}
{\sc G.~Ok{\v s}a, Y.~Yamamoto, M.~Be{\v c}ka, and M.~Vajter{\v s}ic}, {\em
  Asymptotic quadratic convergence of the two-sided serial and parallel
  block-{Jacobi} {SVD} algorithm}, SIAM J. Matrix Anal. Appl., 40 (2019),
  pp.~639--671, \url{https://doi.org/10.1137/18M1222727}.

\bibitem{OksaYamamotoVajtersic2017}
{\sc G.~Ok{\v s}a, Y.~Yamamoto, and M.~Vajter{\v s}ic}, {\em Asymptotic
  quadratic convergence of the serial block-{Jacobi} {EVD} algorithm for
  {H}ermitian matrices}, Numer. Math., 136 (2017), pp.~1071--1095,
  \url{https://doi.org/10.1007/s00211-016-0863-5}.

\bibitem{RoyoEtAl1998}
{\sc D.~Royo, M.~Valero-Garc{\'\i}a, and A.~Gonz{\'a}lez}, {\em A
  {Jacobi}-based algorithm for computing symmetric eigenvalues and eigenvectors
  in a two-dimensional mesh}, in Proceedings of the Sixth Euromicro Workshop on
  Parallel and Distributed Processing (PDP'98), IEEE Computer Society, 1998,
  pp.~463--469, \url{https://doi.org/10.1109/EMPDP.1998.647234}.

\bibitem{Rutishauser1966}
{\sc H.~Rutishauser}, {\em The {Jacobi} method for real symmetric matrices},
  Numer. Math., 9 (1966), pp.~1--10, \url{https://doi.org/10.1007/BF02165223}.

\bibitem{SchwarzEtAl2026}
{\sc A.~Schwarz, R.~M. Larsen, S.~Rodriguez, L.~Kruglyak, and C.~Baker}, {\em
  Mixed-precision {SVD} on {GPUs} via {Ogita--Aishima} iterative refinement},
  ACM Transactions on Mathematical Software,  (2026),
  \url{https://doi.org/10.1145/3860899}.
\newblock Article 3860899, just accepted.

\bibitem{ShroffSchreiber1989}
{\sc G.~Shroff and R.~Schreiber}, {\em On the convergence of the cyclic
  {Jacobi} method for parallel block orderings}, SIAM J. Matrix Anal. Appl., 10
  (1989), pp.~326--346, \url{https://doi.org/10.1137/0610025}.

\bibitem{vanDamHaemersPeek2003}
{\sc E.~R. van Dam, W.~H. Haemers, and M.~B.~M. Peek}, {\em Equitable
  resolvable coverings}, J. Combin. Des., 11 (2003), pp.~113--123,
  \url{https://doi.org/10.1002/jcd.10024}.

\bibitem{ViamonteRalha2006}
{\sc A.~J. Viamonte and R.~Ralha}, {\em High relative precision of eigenvalues
  calculated with {Jacobi} methods}, WSEAS Trans. Math., 5 (2006),
  pp.~721--727.

\bibitem{Wilkinson1962}
{\sc J.~H. Wilkinson}, {\em Note on the quadratic convergence of the cyclic
  {Jacobi} process}, Numer. Math., 4 (1962), pp.~296--300,
  \url{https://doi.org/10.1007/BF01386321}.

\bibitem{ZhouBrent1997}
{\sc B.~B. Zhou and R.~P. Brent}, {\em A parallel ring ordering algorithm for
  efficient one-sided {Jacobi} {SVD} computations}, J. Parallel Distrib.
  Comput., 42 (1997), pp.~1--10, \url{https://doi.org/10.1006/jpdc.1997.1304}.

\bibitem{ZhouBrentKahn1994}
{\sc B.~B. Zhou, R.~P. Brent, and M.~Kahn}, {\em A one-sided {Jacobi} algorithm
  for the symmetric eigenvalue problem}, in Proceedings of the Third Parallel
  Computing Workshop, Kawasaki, Japan, Nov. 1994, pp.~P1--Q--1--P1--Q--7.

\end{thebibliography}

  \FloatBarrier
  \paperbackmatter

\appendix

  \crefalias{section}{appendix}\crefalias{subsection}{subappendix}

\section{Weighted planes}
\label{app:weights}

For a prime power $p<K$ that does not divide $K$, the one-stage filling
of \cref{sec:affine} needs a weighted $K$-arc in $\mathrm{AG}(2,p)$
of large capacity. Let $\mu=\lceil K/p\rceil$ and
$t=\mu p-K$, so that $1\le t<p$. The constant weights $\mu$ give every
line $p\mu=K+t$ columns, $t$ too many. We therefore subtract a
\emph{deficit} $d(P)\ge0$ that removes at least $t$ from every line, that
is, a \emph{$t$-fold blocking multiset}:
$\sum_{P\in\lambda}d(P)\ge t$ for every line $\lambda$. If moreover
$d(P)\le\mu$ at every point, then $w=\mu-d$ is a weighted $K$-arc of
capacity $\mu p^2-\sum_Pd(P)$~\cite[proof of Theorem~3.4]{Furedi1989}.

F\"uredi's multiset has a point of multiplicity $t$~\cite[(1.1)]{Furedi1989},
which exceeds $\mu$ for most orders used here (for example $K=64$, $p=29$,
with $\mu=3$ and $t=23$). We use instead the lines dual to points of a
conic~\cite[Chaps.~7--8]{Hirschfeld1998}, whose multiplicities are at most
two.

\begin{lemma}
\label{lem:weights}
Let $p$ be a prime power, $1\le t<p$, and $\mathcal A\subseteq\F_p$
have $t+1$ elements. For $\alpha\in\mathcal A$, let
$L_\alpha$ be the line $y=\alpha x+\alpha^2$, with arithmetic
in $\F_p$. Then
\[
    d(P)=\bigl|\{\alpha\in\mathcal A:P\in L_\alpha\}\bigr|
\]
is a $t$-fold blocking multiset with total weight $(t+1)p$ and
$0\le d(P)\le2$ at every point.
\end{lemma}

\begin{proof}
At a point $P=(x,y)$, the incident parameters are roots of
$z^2+xz-y$, so there are at most two in any characteristic.
The selected lines have distinct slopes. Any affine line is
parallel to at most one of them and therefore meets at least
$t$ of them. Counting these incidences with multiplicity gives
$\sum_{P\in\lambda}d(P)\ge t$, even when several intersections
coincide. Finally, each of the $t+1$ selected lines contributes
$p$ to the total weight.
\end{proof}

Since $p<K$ gives $\mu\ge2$, this deficit applies for every $t$, and the
weighted $K$-arc has capacity
\begin{equation}
\label{eq:plane-capacity}
    C(p,K)=\mu p^2-(t+1)p=p(K-1).
\end{equation}

  \raggedbottom
  \definecolor{algscope}{HTML}{006FC1}
\definecolor{algsync}{HTML}{E45649}
\definecolor{alghw}{HTML}{009B83}
\newcommand{\Sync}[1]{\State\textcolor{algsync}{\textbf{synchronize #1}}}
\newcommand{\Each}[1]{\textcolor{algscope}{$\forall$\,#1}}
\newcommand{\Hw}[1]{\textcolor{alghw}{#1}}

\section{The inner eigensolver}
\label{sm:evd}
\kgram{} supplies the panel Gram matrix $H_S$ to \kevd{}, which
approximately diagonalizes it in one inner sweep. To check whether
further inner sweeps would save enough outer work to offset their cost,
we tested two and four inner sweeps on five matrices with $n=512$ to $4096$.
In the geometric mean, they reduced the outer sweep count by at
most $5\%$ at $K=64$ and $10\%$ at $K=128$, while taking $1.3$ to
$2.2$ times as long. We therefore fix the inner solve at one sweep
and focus on reducing the cost of that sweep. All three GPU listings
below use $K=64$, with one thread block of eight warps per panel.

\subsection{Blocking the inner sweep}
To identify where a sweep spends its time, start with the scalar
two-sided baseline in \cref{alg:evd-scalar}. It keeps the transformed
Gram matrix $H$ and the accumulated rotations $Z$ in shared memory, a
common choice in batched GPU Jacobi solvers such as that of
MAGMA~\cite{AbdelfattahFasi2026}.
At $K=64$, round-robin visits all column pairs in $63$ stages of
$32$ disjoint pairs, four per warp. Each rotation reads and writes
two columns and two rows of $H$ and two columns of $Z$. Each update takes six
flops for two entries read and written, giving just $3/16$ flop per
byte in double precision. The row updates also depend on entries
just written by the column updates, so each stage requires two block
barriers in this baseline. The sweep thus reads and writes $3K^2$
words per stage and uses $126$ barriers in total.

\begin{algorithm}[tbp]
\small
\caption{Scalar two-sided baseline for \kevd{} at $K=64$.}
\label{alg:evd-scalar}
\begin{algorithmic}[1]
\State \textbf{Input}: $H_S$ ($64\times64$); \textbf{Output}: $Q_S$ ($64\times64$)
\State \textbf{Shared memory}: $H$, $Z$ ($64\times64$)
\State \textbf{Warp registers}: $J_P$ ($2\times2$) per pair
\State \Each{Threads in block} cooperatively: $H\gets H_S$, $Z\gets I_{64}$
\Sync{threads}
\For{stage $k$ : $1\rightarrow63$, with pairs $P_{k,0},\ldots,P_{k,31}$}
  \State \Each{Warp $w$}, for each $j\in\{w,w+8,w+16,w+24\}$:
  \State \quad $P\gets P_{k,j}$, $J_P\gets\textsc{Rotation}(H[P,P])$ \Comment{\Hw{scalar}}
  \State \quad $H[:,P]\gets H[:,P]\,J_P$ \Comment{\Hw{AXPY}}
  \State \quad $Z[:,P]\gets Z[:,P]\,J_P$ \Comment{\Hw{AXPY}}
  \Sync{threads}
  \State \Each{Warp $w$}, for each assigned pair $P$:
  \State \quad $H[P,:]\gets J_P^\top H[P,:]$ \Comment{\Hw{AXPY}}
  \Sync{threads}
\EndFor
\State \Each{Threads in block} cooperatively: $Q_S\gets Z$
\end{algorithmic}
\end{algorithm}

A natural remedy is blocking, as in the outer iteration: a sweep
of two-sided block Jacobi on $H$, with inner pivot sets $L$ of $N$
columns. For each $L$, one warp loads $H[L,L]$ into
registers and uses warp shuffles to perform $N-1$ stages of scalar
Jacobi. It then applies the accumulated rotation $Q_L$ to
$H[:,L]$, $Z[:,L]$ and $H[L,:]$ in three matrix products
(\cref{alg:evd-block}). The products use the FP64 matrix
multiply-accumulate (MMA) instructions of the Hopper tensor cores and
reuse each loaded entry across the $N$ columns or rows of the product. The resulting
arithmetic intensity is $N/8$ flop per byte of shared-memory traffic,
$2N/3$ times that of the scalar updates.

For this blocking to work, the pivot-set size $N$ must meet two
requirements. The $N\times N$ diagonal block must fit in one warp's
registers, and $N$ must be a multiple of the output and reduction
dimensions $8$ and $4$ of the MMA instructions so the products use full tiles. At $K=64$,
$N=8$ meets both requirements and gives a natural schedule based
on $\mathrm{AG}(2,8)$. Assigning one column to each point gives
eight lines of eight columns in each parallel class, one line per
warp. Every pair shares exactly one line, so a sweep takes nine
classes. Block round-robin at the same pivot-set size pairs $16$
fixed blocks of four columns and takes $15$ stages.

Because the eight lines in a class cover all the columns, the
three updates still read and write the same $3K^2$ words as a scalar
stage. The shorter sweep therefore reduces update traffic by a
factor of seven. With $N=8$, the blocked products give one flop
per byte. The shorter sweep also reduces synchronization from
$126$ to $18$ barriers. The two barriers per class still matter
because there is little work between them. The row update for $L$
still reads blocks $H[L,L']$ just written by other warps. Eliminating
that dependency would allow us to use one barrier per class.

\begin{algorithm}[tbp]
\small
\caption{Blocked two-sided baseline for \kevd{} at $K=64$.}
\label{alg:evd-block}
\begin{algorithmic}[1]
\State \textbf{Input}: $H_S$ ($64\times64$); \textbf{Output}: $Q_S$ ($64\times64$)
\State \textbf{Shared memory}: $H$, $Z$ ($64\times64$)
\State \textbf{Warp registers}: $Q_L$ ($8\times8$)
\State \Each{Threads in block} cooperatively: $H\gets H_S$, $Z\gets I_{64}$
\Sync{threads}
\For{class $i$ : $0\rightarrow 8$}
  \State \Each{Warp $w$}: $L\gets$ line $w$ of $\mathcal C_i$
  \State \quad $Q_L\gets\textsc{WarpJacobi}(H[L,L])$ \Comment{\Hw{shuffles}}
  \State \quad $H[:,L]\gets H[:,L]\,Q_L$ \Comment{\Hw{MMA}}
  \State \quad $Z[:,L]\gets Z[:,L]\,Q_L$ \Comment{\Hw{MMA}}
  \Sync{threads}
  \State \Each{Warp $w$}:
  \State \quad $H[L,:]\gets Q_L^\top H[L,:]$ \Comment{\Hw{MMA}}
  \Sync{threads}
\EndFor
\State \Each{Threads in block} cooperatively: $Q_S\gets Z$
\end{algorithmic}
\end{algorithm}

\subsection{One-sided Jacobi on the Gram matrix}
To eliminate the row update, we first consider one-sided Jacobi
on $H_S$~\cite{Hestenes1958}. It stores the accumulated rotations
in $Z$ and applies each rotation to $Y=H_SZ$ from the right.
Each warp therefore updates only its own columns, and a barrier
at the end of a class is sufficient. Since $H_S$ is symmetric
positive semidefinite, its right singular vectors are eigenvectors,
so the converged method diagonalizes $H_S$
\cite{ZhouBrentKahn1994,BeckaOksa2016,BoukaramEtAl2018,OksaBecka2025}.
Our choice of an incomplete solve puts a stronger demand on each
sweep. Its rotations must make enough progress in orthogonalizing
the panel $X_SZ$ to keep the outer sweep count small.

The two methods choose their angles from different matrices. In exact
arithmetic, one-sided Jacobi on a matrix is two-sided Jacobi on its
Gram matrix~\cite[\S4]{DemmelVeselic1992}, here $H_S^2$. With the
current panel Gram matrix $G=Z^\top H_SZ$, two-sided Jacobi chooses the
angle for a pair $\{p,q\}$ to annihilate $g_{pq}$, whereas one-sided
Jacobi on $H_S$ annihilates $(G^2)_{pq}=\sum_k g_{pk}g_{kq}$, which
couples the pair through the other columns. The two angles agree to
first order when $G$ is nearly diagonal with distinct diagonal entries.
While the panel is far from orthogonal, as during a single inner sweep,
one-sided Jacobi can instead rotate a pair that is already orthogonal.

This matters because the outer iteration progresses only as the
off-diagonal part of the Gram matrix shrinks. Let $\operatorname{off}(G)=\sum_{i\neq j}g_{ij}^2$. A rotation
in the plane of $p$ and $q$, applied to the panel columns, changes $G$
to $G'$ with
\[
\operatorname{off}(G')=\operatorname{off}(G)-2g_{pq}^2+2g_{pq}'^2,
\]
since it preserves the Frobenius norm of the $2\times2$ block and of
every pair $(g_{kp},g_{kq})$ with $k\notin\{p,q\}$. The inner products between the panel and
the other columns of $X$ are rotated in the same way, so
$\operatorname{off}(X^\top X)$, which the outer iteration drives to zero,
changes by the same amount. The two-sided angle sets $g_{pq}'=0$, so
every rotation removes $2g_{pq}^2$, the largest decrease possible. The
one-sided angle instead annihilates the $(p,q)$ entry of $G'^2$. It
reduces $\operatorname{off}(G^2)$, but $\operatorname{off}(G)$ only when
$|g_{pq}'|<|g_{pq}|$: when it rotates an orthogonal pair of columns with
different norms, $g_{pq}=0$ and $g_{pq}'\neq0$, and the rotation undoes
progress of the outer iteration.

We measure the effect on the outer sweep count by comparing the two
methods on the host with
$n=128$, $K=32$, and block round-robin. We use two matrices with different
sources of ill-conditioning: a geometric singular spectrum with
$\kappa=10^{10}$ and a graded matrix $D_1BD_2$. Here $B$ has unit
columns and condition number $10$, and $D_1,D_2$ have geometric
diagonal entries with $\kappa_D=10^{10}$ in random order. We also
run one-sided Jacobi to inner convergence, stopping when a sweep
rotates no pair, to see whether further inner work recovers the
lost progress. We measure the inner work by summing the mean
number of inner sweeps per panel over the outer sweeps.

With one inner sweep, one-sided Jacobi on $H_S$ needs $47$ and $75$ outer
sweeps on the two matrices, against $26$ and $25$ for two-sided Jacobi.
Running it to inner convergence recovers the two-sided counts, but at
$8.7$ and $11.0$ times the inner work of one-sweep two-sided Jacobi. \kevd{} therefore needs the
rotations of two-sided Jacobi with the column-only updates of
one-sided Jacobi.

\FloatBarrier

\subsection{One-sided Jacobi for the symmetric eigenproblem}
For the inner solve, we use the block one-sided Jacobi method of
Viamonte and Ralha~\cite{ViamonteRalha2006}, adapted to variable
partitioning with the lines of $\mathrm{AG}(2,8)$ as pivot sets. It
chooses rotations as in two-sided Jacobi and applies them only to
columns. The unblocked method and related one-sided eigensolvers for
symmetric matrices appeared
earlier~\cite{Kaiser1972,Eberlein1987,RoyoEtAl1998,CuencaGimenez2000}.
To our knowledge, the method has yet to be adopted in standard solver
libraries.

The formulation represents the current Gram matrix as $G=Z^\top Y$,
where $Z$ accumulates the rotations and $Y=H_SZ$. To choose
$Q_L$, a warp only needs $B_L=Z[:,L]^\top Y[:,L]=G[L,L]$,
the same diagonal block as in two-sided Jacobi in exact arithmetic.
After choosing $Q_L$ from that block, the warp applies it to the
columns $L$ of both $Y$ and $Z$. This preserves $Y=H_SZ$ for
the next class and confines the updates to the warp's own columns
(\cref{alg:method1}).

\begin{algorithm}[tbp]
\small
\caption{One-sided Jacobi sweep for a symmetric matrix~\cite{ViamonteRalha2006}.}
\label{alg:method1}
\begin{algorithmic}[1]
\State $Z\gets I_K$, $Y\gets H$ \Comment{invariant: $Y=HZ$}
\For{each stage $\mathcal C$ of the sweep}
  \For{each pivot set $L\in\mathcal C$, independently}
    \State $B_L\gets Z[:,L]^\top Y[:,L]$ \Comment{$=(Z^\top HZ)[L,L]$}
    \State $Q_L\gets$ rotations diagonalizing $B_L$ \Comment{$Q_L^\top B_LQ_L\approx\Lambda_L$}
    \State $Z[:,L]\gets Z[:,L]\,Q_L$
    \State $Y[:,L]\gets Y[:,L]\,Q_L$ \Comment{keeps $Y=HZ$}
  \EndFor
\EndFor
\State \Return $Z$
\end{algorithmic}
\end{algorithm}

Because the line updates are independent, \kevd{} keeps the warp
mapping of the blocked sweep. At $K=64$, one thread block visits
the nine classes of $\mathrm{AG}(2,8)$ (\cref{alg:evd-ours}).
Each warp forms its $8\times8$ block $B_L$ with MMA instructions,
computes $Q_L$ by a register-resident Jacobi sweep, and applies it
to $Z[:,L]$ and $Y[:,L]$ in the same way. The last class
can skip the $Y$ update because only $Z$ is returned. One block
barrier per class makes the updated columns visible to the other warps.

\begin{algorithm}[tbp]
\small
\caption{\kevd{} at $K=64$.}
\label{alg:evd-ours}
\begin{algorithmic}[1]
\State \textbf{Input}: $H_S$ ($64\times64$); \textbf{Output}: $Q_S$ ($64\times64$)
\State \textbf{Shared memory}: $Z$, $Y$ ($64\times64$) and per-warp scratch
\State \Each{Threads in block} cooperatively: $Z\gets I_{64}$, $Y\gets H_S$
\Sync{threads}
\For{class $i$ : $0\rightarrow 8$}
  \State \Each{Warp $w$}: $L\gets$ line $w$ of $\mathcal C_i$
  \State \quad $B_L\gets Z[:,L]^\top Y[:,L]$ \Comment{\Hw{MMA}}
  \State \quad $Q_L\gets\textsc{WarpJacobi}\bigl(\tfrac12(B_L+B_L^\top)\bigr)$ \Comment{\Hw{shuffles}}
  \State \quad $Z[:,L]\gets Z[:,L]\,Q_L$ \Comment{\Hw{MMA}}
  \State \quad \textbf{if} $i<8$: $Y[:,L]\gets Y[:,L]\,Q_L$ \Comment{\Hw{MMA}}
  \Sync{threads}
\EndFor
\State \Each{Threads in block} cooperatively: $Q_S\gets Z$
\end{algorithmic}
\end{algorithm}

At $K=128$, the two arrays no longer fit in the shared memory of one
streaming multiprocessor, and each panel is split between two thread
blocks. This variant and further implementation details are documented
in the GitHub repository.\footnote{\coderepository}

\section{Controlling norm drift}
\label{sm:drift}
In the initial version of AG-OSBJ, the backward and singular-value
errors grew with the matrix order $n$, a loss of accuracy similar to
the one Novakovi\'c~\cite[\S5.3.1]{Novakovic2023} reported in his solver.
We identify its source, reduce it during the iteration, and remove what
remains from the returned factors.

\subsection{Diagnosing the accuracy loss}
\label{sm:drift-diagnosis}
On the \mfam{geometric} spectrum with $\kappa=10^{10}$, the backward
error of the initial version was $41$ times that of \texttt{gesvd} at
$n=512$ and about $150$ times at $n=8192$, and its singular-value error
followed the same trend (\cref{tab:postprocessing}). At the same time,
the right factor lost its normalization. Let $W$ be the accumulated
right factor, the array $V$ of \cref{alg:block-jacobi}, which starts
from the identity and accumulates every rotation applied to the
columns. In exact arithmetic, $W$ is orthogonal and each of its columns
$w_j$ has unit norm. We call
\[
\delta_j=\|w_j\|_2^2-1
\]
the \emph{norm drift} of column $j$, and
$\bar\delta=\frac1n\sum_j\delta_j$ the mean drift. In the initial
version, $\bar\delta$ grew roughly in proportion to $n$, from about
$1.2\cdot10^3\,\epsmach$ at $n=512$ to $2.7\cdot10^4\,\epsmach$ at
$n=8192$.

\begin{table}[tbhp]
\caption{Scaled backward error $n\eta$ and singular-value error
$n\rho_S$ of AG-OSBJ with the rotations in the direct form, as in the
initial version, and in Rutishauser's form (\cref{sm:drift-source}):
for the unscaled factors, reconstructed from the returned factors and
the column norms of $W$, and for the returned, rescaled ones
(\cref{sm:drift-rescaling}); consistency error $\|F\|_F/\|A\|_F$;
\texttt{gesvd} for reference. Error measures as in
\cref{sec:exp-cusolver}; \mfam{geometric} spectrum, $\kappa=10^{10}$,
$K=64$.}
\label{tab:postprocessing}
\centering
\small
\begin{tabular}{@{}lccccc@{}}
\toprule
& \multicolumn{2}{c}{Unscaled} & \multicolumn{2}{c}{Returned} & \\
\cmidrule(lr){2-3} \cmidrule(lr){4-5}
Method & $n\eta$ & $n\rho_S$ & $n\eta$ & $n\rho_S$ & $\|F\|_F/\|A\|_F$ \\
\midrule
\multicolumn{6}{@{}l}{$n=512$} \\
\quad\texttt{gesvd} & -- & -- & $6.8\cdot10^{-15}$ & $3.0\cdot10^{-15}$ & -- \\
\quad Direct form & $2.8\cdot10^{-13}$ & $1.4\cdot10^{-13}$ & $8.9\cdot10^{-15}$ & $7.4\cdot10^{-16}$ & $6.0\cdot10^{-15}$ \\
\quad Rutishauser's form & $2.4\cdot10^{-14}$ & $1.1\cdot10^{-14}$ & $6.8\cdot10^{-15}$ & $3.7\cdot10^{-16}$ & $6.1\cdot10^{-15}$ \\
\addlinespace[4pt]
\multicolumn{6}{@{}l}{$n=8192$} \\
\quad\texttt{gesvd} & -- & -- & $3.2\cdot10^{-14}$ & $1.8\cdot10^{-14}$ & -- \\
\quad Direct form & $4.9\cdot10^{-12}$ & $2.5\cdot10^{-12}$ & $6.6\cdot10^{-14}$ & $1.7\cdot10^{-14}$ & $3.2\cdot10^{-14}$ \\
\quad Rutishauser's form & $4.4\cdot10^{-13}$ & $2.2\cdot10^{-13}$ & $2.9\cdot10^{-14}$ & $7.3\cdot10^{-16}$ & $2.5\cdot10^{-14}$ \\
\bottomrule
\end{tabular}

\end{table}

To see how the backward error and the drift connect, we split the
residual into parts that can be measured separately. Let $X$ be the
working matrix, so that the unscaled factors satisfy $U_0\Sigma_0=X$
and $V_0=W$. We define the \emph{consistency error} $F=X-AW$ and the
loss of orthogonality $\Delta=W^\top W-I$, whose diagonal entries are
the drifts $\delta_j$. We also write $\omega_j=\|w_j\|_2$, so that
$\delta_j=\omega_j^2-1$.

\begin{lemma}
\label{lem:residual}
If $\|\Delta\|_2<1$, then $W$ is nonsingular,
$\|W^{-1}\|_2\le(1-\|\Delta\|_2)^{-1/2}$, and
\begin{equation}
\label{eq:backward-identity}
A-U_0\Sigma_0V_0^\top
=-\bigl(X\operatorname{diag}(\delta_j)
+X(\Delta-\operatorname{diag}(\delta_j))+F\bigr)W^{-1}.
\end{equation}
\end{lemma}
\begin{proof}
The residual is $A-XW^\top$, and
$(A-XW^\top)W=AW-X(I+\Delta)=-F-X\Delta$. The squared singular values of
$W$ are the eigenvalues of $I+\Delta$.
\end{proof}

The identity separates three sources of error: the norm drift, the
off-diagonal loss of orthogonality, and the consistency error. In every
run of this section, $\|\Delta\|_F<10^{-9}$, so the factor $W^{-1}$ is
negligible, and the consistency error stays within $12\%$ of the
backward error of \texttt{gesvd}. The drift term
$\|\Sigma_0\operatorname{diag}(\delta_j)\|_F/\|A\|_F$, by contrast,
reproduces the backward error of the initial version to within $1\%$
from the fourth sweep on (\cref{fig:drift}). The accuracy loss
therefore comes from norm drift: the columns of $W$, and with them
those of $X$, acquire incorrect lengths.

\subsection{Reducing norm drift}
\label{sm:drift-source}
To see where the drift comes from, consider one stage that applies a
computed transformation $\hat Q$ to both arrays: $X^+=X\hat Q+E_X$ and
$W^+=W\hat Q+E_W$, where $E_X$ and $E_W$ are the rounding errors of the
two products in \krot{}. Let $C=FW^{-1}=XW^{-1}-A$ denote the
consistency contribution to \eqref{eq:backward-identity}.

\begin{lemma}
\label{lem:stage}
If $W$ and $W^+$ are nonsingular, then
\begin{align}
\label{eq:consistency}
C^+&=C+\bigl(E_X-(A+C)E_W\bigr)(W^+)^{-1},\\
\label{eq:drift}
\Delta^+&=\hat Q^\top\Delta\hat Q+(\hat Q^\top\hat Q-I)
+(W\hat Q)^\top E_W+E_W^\top W\hat Q+E_W^\top E_W.
\end{align}
\end{lemma}
\begin{proof}
Since $X=(A+C)W$, we have $X^+=(A+C)(W^+-E_W)+E_X$, which gives
\eqref{eq:consistency}. Expanding $(W\hat Q+E_W)^\top(W\hat Q+E_W)-I$
with $W^\top W=I+\Delta$ gives \eqref{eq:drift}.
\end{proof}

The nonorthogonality of $\hat Q$ enters $\Delta$ directly, even with
exact matrix products, whereas $C$ changes only through their rounding
errors. Computing the rotations more accurately can therefore reduce
the drift without affecting the consistency error. To see how the
rotations contribute, we relate the mean drift to their coefficients.
Transformations accumulate at three levels, each starting
from the identity: the warp transformation $Q_L$ accumulates plane
rotations, the panel transformation $Q_S$ accumulates the $Q_L$, and
the right factor $W$ accumulates the stage transformations. For a
rounded update $Y^+=YT+E$ of a square array, expanding the squared
Frobenius norm gives
\begin{align}
\|Y^+\|_F^2-\|Y\|_F^2
={}&\operatorname{tr}(T^\top T-I)+2\langle YT,E\rangle \notag\\
&+\operatorname{tr}\bigl((Y^\top Y-I)(TT^\top-I)\bigr)+\|E\|_F^2,
\label{eq:trace-step}
\end{align}
where $\langle\cdot,\cdot\rangle$ is the Frobenius inner product.
For a left update $Y^+=TY+E$, the same identity holds with
$TY$, $YY^\top$, and $T^\top T$ replacing $YT$, $Y^\top Y$,
and $TT^\top$ in the last three terms.

Applying this identity successively at the three levels expresses
the mean drift $\bar\delta$ in terms of
the defects of the individual plane rotations. Summing over
the updates of $W$ gives $\|W\|_F^2-n=n\bar\delta$, while the
trace contribution of each stage is the sum of $\|Q_S\|_F^2-K$
over its panels. Expanding these contributions through the updates
of $Q_S$, and then those of $Q_L$, leaves the trace defects of the
individual rotations. For a plane rotation with coefficients $a,b$,
this defect is $2d$, where $d=a^2+b^2-1$. Consequently,
\begin{equation}
\label{eq:trace}
n\bar\delta=2\sum d+r,
\end{equation}
where the sum runs over all rotations applied so far. This identity
is exact: the remainder $r$ collects the rounding-error terms,
products of losses of orthogonality, and squared rounding errors
from \eqref{eq:trace-step} at all three levels.

\begin{figure}[tbhp]
\centering
\definecolor{ours}{HTML}{1B4F8A}
\definecolor{robin}{HTML}{C8553D}
\definecolor{bound}{HTML}{5A5A5A}
\definecolor{gesvdj}{HTML}{3A7D44}
\definecolor{gesvdp}{HTML}{C99A2E}
\begin{tikzpicture}
\begin{groupplot}[
    group style={group size=2 by 1, horizontal sep=1.6cm},
    width=0.5\linewidth, height=0.42\linewidth,
    scale only axis=false,
    axis lines=left, axis line style={bound!60, line width=1bp, line cap=round},
    tick style={bound!60, line width=1bp}, tick align=outside,
    major tick length=2.5pt,
    ticklabel style={font=\footnotesize, text=bound},
    /pgf/number format/1000 sep={},
    label style={font=\small},
    every axis title/.style={font=\small, at={(0.5,1.04)}, anchor=south},
    grid=major, major grid style={bound!15, line width=1bp},
    every axis plot/.append style={line join=round},
    legend cell align=left,
    legend style={draw=none, fill=none,
        font=\scriptsize, row sep=-2pt, inner sep=1pt},
]
\nextgroupplot[
    title={Direct form},
    xmin=0,
    xlabel={Sweep},
    ymode=log,
    ymin=1e-16,
    ymax=1e-11,
    ytick={1e-16,1e-15,1e-14,1e-13,1e-12,1e-11},
    legend pos=south east,
]
\addplot[ours, solid, line width=1.1bp] table[x=sweep, y=drift] {figures/drift_direct.dat};
\addlegendentry{Mean norm drift}
\addplot[robin, only marks, mark=o, mark size=1.2pt, line width=1bp] table[x=sweep, y=predicted] {figures/drift_direct.dat};
\addlegendentry{Coefficient contribution}
\addplot[bound, dashed, line width=1.1bp] table[x=sweep, y=consistency] {figures/drift_direct.dat};
\addlegendentry{Consistency error}
\nextgroupplot[
    title={Rutishauser's form},
    xmin=0,
    xlabel={Sweep},
    ymode=log,
    ymin=1e-16,
    ymax=1e-11,
    ytick={1e-16,1e-15,1e-14,1e-13,1e-12,1e-11},
    legend pos=north west,
    yticklabels={},
]
\addplot[ours, solid, line width=1.1bp] table[x=sweep, y=drift] {figures/drift_rutishauser.dat};
\addlegendentry{Mean norm drift}
\addplot[bound, dashed, line width=1.1bp] table[x=sweep, y=consistency] {figures/drift_rutishauser.dat};
\addlegendentry{Consistency error}
\end{groupplot}
\end{tikzpicture}
\caption{Mean norm drift and consistency error after each sweep of one
solve, in the direct form (left) and Rutishauser's form (right); on the
left also the coefficient contribution $2\sum d_D/n$ to the drift.
\mfam{Geometric} spectrum, $\kappa=10^{10}$, $n=4096$, $K=64$.}
\label{fig:drift}
\end{figure}
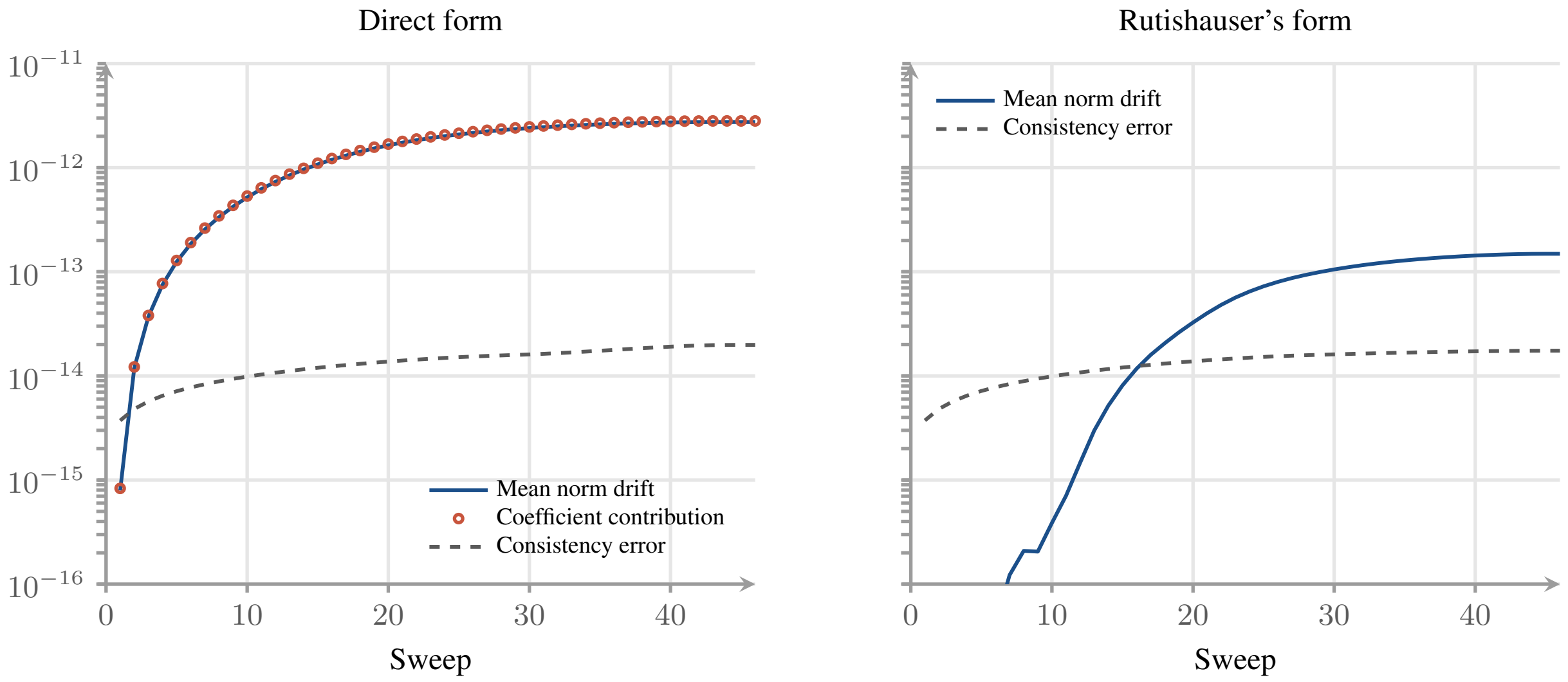

The initial version applied each rotation of the pivot columns $X_p$
and $X_q$, with computed cosine $\hat c$ and sine $\hat s$, in the
standard \emph{direct form}~\cite{GolubVanLoan2013},
\[
X_p'=\hat c X_p-\hat s X_q,\qquad
X_q'=\hat s X_p+\hat c X_q,
\]
whose coefficient defect is $d_D=\hat c^2+\hat s^2-1$. Recorded in
double-double arithmetic, the coefficient contribution $2\sum d_D/n$
reproduces the observed mean drift to within $5\%$ in every sweep of
\cref{fig:drift}. The defects are largest for small rotations,
averaging about $\epsmach/4$ for $10^{-8}<|\hat s|\le10^{-4}$.

To suppress these defects, we apply the rotations in
\emph{Rutishauser's form}~\cite{Rutishauser1966,Drmac1997}, which is
equivalent in exact arithmetic but writes the rotation as a correction
to the identity:
\[
X_p'=X_p-\hat s(X_q+\hat\tau X_p),\qquad
X_q'=X_q+\hat s(X_p-\hat\tau X_q),\qquad
\hat\tau=\operatorname{fl}\!\left(\frac{\hat s}{1+\hat c}\right).
\]
Its coefficient defect is $d_R=(1-\hat s\hat\tau)^2+\hat s^2-1$. With
the exact quotient $\tau_0=\hat s/(1+\hat c)$, this becomes
$d_R=\tau_0^2d_D$. Since $\hat c\ge0$, we have $|\tau_0|\le|\hat s|$,
and a relative error of order $\epsmach$ in $\hat\tau$ perturbs
$1-\hat s\hat\tau$ by $O(\hat s^2\epsmach)$. A direct-form defect of
order $\epsmach$ thus becomes $O(\hat s^2\epsmach)$.

In the run of \cref{fig:drift}, Rutishauser's form reduces the mean
drift by a factor of eighteen, to $670\,\epsmach$, while the consistency
error barely changes: across the five sizes, it differs from that of the
direct form by at most $22\%$. Its coefficient contribution
$2\sum d_R/n$ is only $-0.1\,\epsmach$, so the remaining drift comes
from the other rounding errors in $r$. Rather than attempting to
eliminate each of them, we remove the remaining drift from the returned
factors.

\subsection{Normalizing the final factors}
\label{sm:drift-rescaling}
We normalize each column $w_j$ by its length $\omega_j$ and apply
the same scaling to $x_j$. Since
$x_j/\omega_j=A(w_j/\omega_j)+f_j/\omega_j$, this preserves the
consistency relation up to the same scaling. The resulting factors are
\begin{equation}
\label{eq:rescaling}
\hat v_j=\frac{w_j}{\omega_j},\qquad
\hat\sigma_j=\frac{\|x_j\|_2}{\omega_j},\qquad
\hat u_j=\frac{x_j}{\|x_j\|_2}.
\end{equation}
The rescaling requires one additional pass over $W$.\footnote{LAPACK's
\texttt{DGESVJ} also normalizes its right vectors~\cite{DrmacVeselic2008b}.}

To see how this changes the backward error, consider exact rescaling
of the stored arrays. With $\Omega=\operatorname{diag}(\omega_j)$,
we have $\hat V=W\Omega^{-1}$ and $\hat U\hat\Sigma=X\Omega^{-1}$.
Using $\omega_j^2=1+\delta_j$ in the new loss of orthogonality $E_V=\hat V^\top\hat V-I$
and applying the residual identity \eqref{eq:backward-identity} gives
\begin{gather*}
E_V=\Omega^{-1}\bigl(\Delta-\operatorname{diag}(\delta_j)\bigr)\Omega^{-1},\\
A-\hat U\hat\Sigma\hat V^\top
=-(\hat U\hat\Sigma E_V+F\Omega^{-1})\hat V^{-1}.
\end{gather*}
The zero diagonal of $E_V$ removes the drift term from the residual.
The off-diagonal loss of orthogonality remains, and the consistency contribution
$F\Omega^{-1}\hat V^{-1}=FW^{-1}$ is exactly preserved. The left
factor stays $U_0$ because positive rescaling preserves the
direction of each $x_j$.

The effect on singular values is clearest when the error is purely a
stretch. In this case, $W=V_*\Omega$ and $X=U_*\Sigma_*\Omega$ for an
exact SVD $A=U_*\Sigma_*V_*^\top$. An unscaled singular-value estimate
$\omega_j\sigma_j$ then has relative error
$\omega_j-1=\delta_j/2+O(\delta_j^2)$. In the reconstruction,
both factors carry the stretch, giving a relative error
$\omega_j^2-1=\delta_j$ in the corresponding singular direction.
This explains why the normwise singular-value error is about half
the backward error in \cref{tab:postprocessing}. Dividing out
$\omega_j$ recovers the exact factors in the pure-stretch case.

With Rutishauser's form, rescaling reduces the singular-value error by
factors of $31$ at $n=512$ and $302$ at $n=8192$, and brings the
backward error to within $15\%$ of the consistency error, close to that
of \texttt{gesvd} (\cref{tab:postprocessing}).

Rescaling alone would not suffice, because later rotations move drift
off the diagonal, where it survives rescaling. For an orthogonal
rotation of columns $p$ and $q$ with cosine $c$ and sine $s$,
\[
(Q^\top\Delta Q)_{pq}=(\delta_p-\delta_q)cs+(c^2-s^2)\Delta_{pq},
\]
so unequal drifts create an off-diagonal loss of orthogonality even
when $\Delta_{pq}=0$. At $n=8192$, the normalized direct form leaves
seven times the loss of orthogonality $\|E_V\|_F$ of Rutishauser's
form, and its backward error is about twice its consistency error.

\FloatBarrier

\end{document}